\documentclass[a4paper,11pt]{amsart}
\usepackage{amssymb,amsfonts,amsmath}
\def\d{\delta}

\def\C{\mathbb{C}}
\def\c2{\mathbb{C}^2}
\def\R{\mathbb{R}}

\def\Z{\mathbb{Z}}
\def\N{\mathbb{N}}
\def\P{\mathbb{P}}

\def\1{\bold{1}}

\def\a{\alpha}

\def\e{\varepsilon}
\def\l{\lambda}
\def\f{\varphi}
\def\g{\gamma}

\def\p{\psi}
\def\r{\varrho}
\def\om{\omega}

\def\E{{\mathcal{E}_{\chi}(X,\om)}}

\newtheorem{lem}{Lemma}[section]
\newtheorem{pro}[lem]{Proposition}
\newtheorem{defi}[lem]{Definition}
\newtheorem{def/not}[lem]{Definition/Notations}

\newtheorem{thm}[lem]{Theorem}
\newtheorem{cor}[lem]{Corollary}
\newtheorem{rqe}[lem]{Remark}

\newtheorem{exa}[lem]{Example}

\newenvironment{proof3.1}
{\noindent {\it{Proof of theorem 3.1}}}{$\Box$ \linebreak[4]}

\begin{document}

\title[The Weighted Monge-Amp\`ere energy of qpsh functions]
{The Weighted Monge-Amp\`ere energy of quasiplurisubharmonic functions}

\author{Vincent GUEDJ \& Ahmed ZERIAHI}

\begin{abstract}
We study degenerate complex Monge-Amp\`ere equations on 
a compact K\"ahler manifold $(X,\om)$.
We show that the complex Monge-Amp\`ere operator
$(\om+dd^c \cdot)^n$ is well-defined on the class
${\mathcal E}(X,\om)$ of $\om$-plurisubharmonic functions
with finite weighted Monge-Amp\`ere energy.
The class ${\mathcal E}(X,\om)$ is the largest class of 
$\om$-psh functions on which the Monge-Amp\`ere operator is 
well-defined and the comparison principle is valid. It
contains several functions whose gradient is
not square integrable. We give a complete description of 
the range of the operator $(\om+dd^c \cdot)^n$
on ${\mathcal E}(X,\om)$, as well as on some of its subclasses.

We also study uniqueness properties, extending Calabi's result to this
unbounded and degenerate situation, and
we give applications to complex dynamics and to the existence of 
singular K\"ahler-Einstein
metrics.
\end{abstract}

\maketitle


\section*{Introduction}

Let $X$ be a compact connected K\"ahler manifold of complex dimension $n \in \N^*$.
Let $\om$ be a K\"ahler form on $X$.
Given $\mu$ a positive Radon  measure on $X$ such that $\mu(X)=\int_X \om^n$,
we study the complex Monge-Amp\`ere
equation
$$
(MA)_\mu \hskip2cm (\om+dd^c \f)^n=\mu,
$$
where $\f$, the unknown function, is such that $\om_{\f}:=\om+dd^c \f$ is a positive current.
Such functions are called $\om$-plurisubharmonic
($\om$-psh for short).
We refer the reader to [GZ 1] for 
basic properties of the set $PSH(X,\om)$ of all such functions.
Here $d=\partial+\overline{\partial}$
and $d^c=\frac{1}{2 i \pi} (\partial-\overline{\partial})$.

Complex Monge-Amp\`ere equations have been studied by several authors over the 
last fifty years, in connection with questions from K\"ahler geometry and
complex dynamics (see [A], [Y], [T], [K 1,2], [DP], [S], [EGZ] for references).
The first and cornerstone result is due to S.T.Yau who proved [Y]
that $(MA)_\mu$ admits a solution $\f \in PSH(X,\om) \cap {\mathcal C}^{\infty}(X)$
when $\mu=f \om^n$ is a smooth volume form.

Motivated by applications towards complex dynamics, we need here to 
consider measures $\mu$ which are quite singular, whence to deal with singular
$\om$-psh functions $\f$.
We introduce and study a class ${\mathcal E}(X,\om)$ of $\om$-psh functions
for which the complex Monge-Amp\`ere operator
$(\om+dd^c \f)^n$ is well-defined (see Definition 1.1): 
following E.Bedford and A.Taylor [BT 4] we show
that the operator $(\om+dd^c \f)^n$ is well defined in
$X \setminus (\f=-\infty)$ for all functions
$\f \in PSH(X,\om)$; the class ${\mathcal E}(X,\om)$ is the set of
functions $\f \in PSH(X,\om)$ such that $(\om+dd^c \f)^n$
has full mass $\int_X \om^n$ in $X \setminus (\f=-\infty)$.
When $n=\dim_{\C} X=1$, this is precisely the subclass
of functions $\f \in PSH(X,\om)$ whose Laplacian does not charge
polar sets. It is striking that the class
${\mathcal E}(X,\om)$ contains many functions whose
gradient is not square integrable, hence several results to follow
have no local analogue (compare [Bl 2,3]).

One of our main results gives a complete
characterization of the range of the complex Monge-Amp\`ere
operator on the class ${\mathcal E}(X,\om)$.
\vskip.2cm

\noindent {\bf Theorem A.} 
{\it There exists $\f\in {\mathcal E}(X,\om)$ such that $\mu=(\om+dd^c \f)^n$
if and only if $\mu$ does not charge pluripolar sets.}
\vskip.2cm

An important tool we use is the {\it comparison principle} that we establish
in section 1: we show that ${\mathcal E}(X,\om)$ is the largest class
of $\om$-psh functions on which the complex Monge-Amp\`ere
operator $(\om+dd^c \cdot)^n$ is well defined and the comparison principle
is valid. Another crucial tool for our study is the notion of 
{\it weighted Monge-Amp\`ere energy},
defined as 
$$
E_{\chi}(\f):=\int_X (-\chi) \circ \f \, (\om+dd^c \f)^n,
$$
where $\chi:\R^- \rightarrow \R^-$ is an increasing function such that
$\chi(-\infty)=-\infty$. The properties of this energy are quite different
whether the weight $\chi$ is convex ($\chi \in {\mathcal W}^-$)
or concave ($\chi \in {\mathcal W}^+$). We show (Proposition 2.2) that
$$
{\mathcal E}(X,\om)=\bigcup_{\chi \in {\mathcal W}^- } {\mathcal E}_{\chi}(X,\om),
$$
where ${\mathcal E}_{\chi}(X,\om)$ denotes the class of functions
$\f \in {\mathcal E}(X,\om)$ such that 
$\chi (\f-\sup_X \f) \in L^1((\om+dd^c \f)^n)$.
At the other extreme, we show (Proposition 3.1) that
$$
PSH(X,\om) \cap L^{\infty}(X)=\bigcap_{\chi \in {\mathcal W}^+} {\mathcal E}_{\chi}(X,\om).
$$

The function $\chi(t)=t$ is -- up to multiplicative constant -- the unique weight in
${\mathcal W}^- \cap {\mathcal W}^+$. We let
$$
{\mathcal E}^1(X,\om):=\{ \f \in {\mathcal E}(X,\om) \,  / \, 
\f \in L^1( (\om+dd^c \f)^n) \}
$$
denote the class ${\mathcal E}_{\chi}(X,\om)$ for $\chi(t)=t$.
When $n=\dim_{\C} X=1$, this is the classical class of quasi-subharmonic functions
of finite energy. It deserves special attention both for the theory and the 
applications. Indeed all functions $\f \in {\mathcal E}^1(X,\om)$ have gradient
in $L^2(\om^n)$, while the gradient of  most functions in
${\mathcal E}_{\chi}(X,\om)$, $\chi \in {\mathcal W}^- \setminus {\mathcal W}^+$,
does not belong to $L^2(\om^n)$ (see Example 2.14), in contrast
with the local theory [Bl 2,3].

We obtain the following extension of Calabi's uniqueness result [Ca].

\vskip.2cm
\noindent {\bf Theorem B.}
{\it Assume 
$(\om+dd^c \f)^n \equiv (\om+dd^c \p)^n$ with $\f \in {\mathcal E}(X,\om)$
and $\p \in {\mathcal E}^1(X,\om)$. 
Then $\f-\p$ is constant.}
\vskip.2cm

It is an interesting open question to prove uniqueness of solutions
as above in the larger class ${\mathcal E}(X,\om)$.

We also study the range of the Monge-Amp\`ere operator 
on subclasses ${\mathcal E}_{\chi}(X,\om)$, when
$\chi(t)=-(-t)^p$, $p>0$, letting
$$
{\mathcal E}^p(X,\om):=\{ \f \in {\mathcal E}(X,\om) \, / \,
\f \in L^p((\om+dd^c \f)^n) \}
$$
denote the corresponding class ${\mathcal E}_{\chi}(X,\om)$.

\vskip.2cm
\noindent {\bf Theorem C.}
{\it There exists $\f \in {\mathcal E}^p(X,\om)$ such that 
$\mu=(\om+dd^c \f)^n$ if and only if ${\mathcal E}^p(X,\om) \subset L^p(\mu)$.}
\vskip.2cm

This result also holds for quasihomogeneous weights (see Lemma 3.9).
It was obtained in a local context by U.Cegrell [Ce 1] when $p \geq 1$.

On our way to prove Theorems A,B,C, we establish several intermediate
results (the class ${\mathcal E}(X,\om)$ is convex, weighted Monge-Amp\`ere
operators are continuous under decreasing sequences, etc)
most of which are valid when $\om$ is merely a positive closed $(1,1)$-current
with bounded potentials.
This is motivated by applications towards complex dynamics and complex
differential geometry, as we briefly indicate in section 5.
\vskip.2cm

This article is an expanded version of our previous preprint [GZ 2].

\noindent {\bf Acknowledgement.} We are grateful to the referee 
for his careful reading and for his suggestions
which helped to clarify the exposition.

\section{The class ${\mathcal E}(X,\om)$}

In the sequel $X$ is a compact K\"ahler manifold of dimension $n$,
and $\om$ is a positive closed $(1,1)$-current with bounded potentials, 
such that $\int_X \om^n >0$.
We let $PSH(X,\om)=\{ \f \in L^1(X) \, / \, \f \text{ is u.s.c. and } dd^c \f \geq -\om \}$
denote the set of $\om$-plurisubharmonic functions ($\om$-psh for short)
which was introduced and studied in [GZ 1].

\subsection{Defining the complex Monge-Amp\`ere operator}

It follows from their plurifine properties that if $u,v$ are bounded
plurisubharmonic functions in some open subset $D$ of $\C^n$, then
\begin{equation}
\1_{\{u>v\}} [dd^c u]^n=\1_{\{u>v\}} [dd^c \max (u,v)]^n
\end{equation}
in the sense of Borel measures in $D$ (see [BT 4]).

Let $\f$ be some unbounded $\om$-psh function on $X$ and consider
$\f_j:=\max(\f,-j) \in PSH(X,\om)$ the {\it canonical approximation} of $\f$ by
bounded $\om$-psh functions. This is a decreasing sequence such that, by (1),
$$
\1_{\{\f_j>-k\}} [\om+dd^c \f_j]^n=\1_{\{\f_j>-k\}} [\om+dd^c \max(\f_j,-k)]^n.
$$
Now if $j \geq k$, then $(\f_j>-k)=(\f>-k)$
and $\max(\f_j,-k)=\f_k$, thus
\begin{equation}
\1_{\{\f>-k\}} [\om+dd^c \f_j]^n=\1_{\{\f>-k\}} [\om+dd^c \f_k]^n.
\end{equation}
Observe also that $(\f>-k) \subset (\f>-j)$, therefore
$$
j \geq k \Longrightarrow 
\1_{\{\f>-j\}} [\om+dd^c \f_j]^n \geq \1_{\{\f>-k\}} [\om+dd^c \f_k]^n,
$$
in the weak sense of Borel measures.
Since the total mass of the measures $\1_{\{\f>-j\}} [\om+dd^c \f_j]^n$
is uniformly bounded from above by $\int_X \om^n$,
by Stokes theorem, we can define
$$
\mu_{\f}:=\lim_{j \rightarrow +\infty} \1_{\{\f>-j\}} [\om+dd^c \f_j]^n.
$$

This is a positive Borel measure which is precisely the non-pluripolar part
of $(\om+dd^c \f)^n$, as considered in a local context
by E.Bedford and A.Taylor in [BT 4]. 
Its total mass $\mu_{\f}(X)$ can take any value in $[0,\int_X \om^n]$.

\begin{defi}
We set
$$
{\mathcal E}(X,\om):=\{ \f \in PSH(X,\om) \, / \, \mu_{\f}(X)=\int_X \om^n \}.
$$
\end{defi}

An alternative definition is given by the following observation:

\begin{lem}
Let $(s_j)$ be any sequence of real numbers converging to $+\infty$, such that
$s_j \leq j$ for all $j \in \N$.
The following conditions are equivalent:

(a)  $\f \in {\mathcal E}(X,\om)$;

(b) $(\om+dd^c \f_j)^n(\f \leq -j) \rightarrow 0$;

(c) $(\om+dd^c \f_j)^n(\f \leq -s_j) \rightarrow 0$,

\noindent where $\f_j:=\max(\f,-j)$ denotes the canonical approximation.
\end{lem}

\begin{proof}
By definition of $\mu_{\f}$, we have
$$
\mu_{\f}(X)=\int_X \om^n-\lim_{j \rightarrow +\infty} (\om+dd^c \f_j)^n(\f \leq -j),
$$
hence (a) is equivalent to (b).
It follows from (2) (with $k=s_j$) that
$$
(\om+dd^c \f_j)^n(\f \leq -s_j)=(\om+dd^c \f_{s_j})^n(\f \leq -s_j),
$$
since $s_j \leq j$.
This shows that $(b)$ is equivalent to (c).
\end{proof}

When $X$ is a compact Riemann surface ($n=\dim_{\C} X=1$),
the set ${\mathcal E}(X,\om)$ is the set of $\om$-subharmonic functions
whose Laplacian does not charge polar sets.
It follows from the above discussion that any $\om$-psh function has
a well defined complex Monge-Amp\`ere operator 
$(\om+dd^c \f)^n$ in $X \setminus (\f=-\infty)$, and $\mu_{\f}$ is the trivial
extension of $(\om+dd^c \f)^n$ through $(\f=-\infty)$.
A function $\f$ belongs to ${\mathcal E}(X,\om)$ precisely when
its complex Monge-Amp\`ere has total mass $\int_X \om^n$ in 
$X \setminus (\f=-\infty)$, hence it is natural to use
the notation
$$
(\om+dd^c \f)^n:=\mu_{\f}=\lim_{j \rightarrow +\infty} \1_{\{\f>-j\}} [\om+dd^c \f_j]^n,
$$
for $\f \in {\mathcal E}(X,\om)$.

\begin{thm}
Let $\f \in {\mathcal E}(X,\om)$. Then for all bounded Borel function $b$,
\begin{equation}
\langle (\om+dd^c \f)^n,b \rangle=\lim_{j \rightarrow +\infty}
\langle (\om+dd^c \f_j)^n,b \rangle,
\end{equation}
where $\f_j:=\max(\f,-j)$ is the canonical approximation of $\f$.

In particular $(\om+dd^c \f)^n$ puts no mass on pluripolar sets, and
\begin{equation}
\1_{\{\f>-j\}} (\om+dd^c \f)^n(B)=\1_{\{\f>-j\}} (\om+dd^c \f_j)^n(B)
\end{equation}
for all Borel subsets $B \subset X$.
\end{thm}

Let us emphasize that the convergence in (3) implies -- but is much stronger than --
the convergence in the weak sense of positive Borel measures,
$$
(\om+dd^c \f_j)^n \longrightarrow (\om+dd^c \f)^n.
$$
Several results to follow are a consequence of fact (3): the complex
Monge-Amp\`ere measure $(\om+dd^c \f)^n$ of a function
$\f \in {\mathcal E}(X,\om)$ is very well approximated by the
Monge-Amp\`ere measures $(\om+dd^c \f_j)^n$.

Note also that the complex Monge-Amp\`ere operator 
$(\om+dd^c \f)^n$ is thus well defined for functions
$\f \in {\mathcal E}(X,\om)$, although these functions need not
have gradient in $L^2(X)$ (see Theorem 1.9 and Example 2.14). 
This is in contrast with
the local theory [Bl 2,3].

\begin{proof}
Recall that
$\f \in {\mathcal E}(X,\om)$ if and only if
$(\om+dd^c \f_j)^n (\f \leq -j) \rightarrow 0$,
where $\f_j:=\max(\f,-j)$.
We infer that for all Borel subset $B \subset X$,
$$
(\om+dd^c \f)^n(B):=\lim_{j \rightarrow +\infty} \int_{B \cap (\f >-j)} (\om+dd^c \f_j)^n
=\lim_{j \rightarrow +\infty} \int_{B} (\om+dd^c \f_j)^n.
$$
This yields (3), by using Lebesgue dominated convergence theorem.

Since $(\om+dd^c \f_j)^n$ does not charge any pluripolar set $B$, the same
property holds for $(\om+dd^c \f)^n$.
Equality (4) now follows from (2) and (3).
\end{proof}

Since $(\om+dd^c \f)^n$ does not charge the pluripolar set $(\f=-\infty)$
when $\f \in {\mathcal E}(X,\om)$, one can 
construct a continuous increasing function
$h:\R^+ \rightarrow \R^+$ such that $h(+\infty)=+\infty$ and
$h \circ |\f| \in L^1((\om+dd^c \f)^n)$.
This motivates the following result.

\begin{pro}
Fix $\f \in {\mathcal E}(X,\om)$ and
let $h:\R^+ \rightarrow \R^+$ be a continuous increasing function such that
$h(+\infty)=+\infty$. Then
$$
\int_X h \circ |\f| (\om+dd^c \f)^n <+\infty \Longleftrightarrow 
\sup_{j \geq 0} \int_X h \circ |\f_j| (\om+dd^c \f_j)^n <+\infty,
$$
where $\f_j:=\max(\f,-j)$. 

Moreover if this condition holds, then for all Borel subset $B \subset X$,
\begin{equation}
\int_B h \circ |\f_j| (\om+dd^c \f_j)^n \longrightarrow 
\int_B h \circ |\f| (\om+dd^c \f)^n.
\end{equation}
\end{pro}

\begin{proof}
We can assume without loss of generality that $\f,\f_j \leq 0$.

Assume first that $\sup_{j \geq 0} \int_X h \circ |\f_j| (\om+dd^c \f_j)^n <+\infty$.
Since the Borel measures $h(-\f_j) (\om+dd^c \f_j)^n$ have uniformly bounded masses,
they form a weakly compact sequence. Let $\nu$ be a cluster point.
Since the functions $h (-\f_j)$ increase towards $h(-\f)$
and $(\om+dd^c \f_j)^n$ converges towards $(\om+dd^c \f)^n$, it follows
from semi-continuity that
$h(-\f) (\om+dd^c \f)^n \leq \nu$, hence
$\int_X h(-\f) (\om+dd^c \f)^n \leq \nu(X)<+\infty$.

Conversely assume that $h(-\f) \in L^1((\om+dd^c \f)^n)$. It follows from (4) that 
$$
\int_{(\f \leq -j)} \om_{\f_j}^n=\int_X \om_{\f_j}^n-\int_{(\f>-j)} \om_{\f_j}^n
=\int_{(\f \leq-j)} \om_{\f}^n.
$$
Here -- and in the sequel -- we use the notation $\om_{\f}:=\om+dd^c \f$,
$\om_{\f_j}=\om+dd^c \f_j$. Thus by (4) again,
\begin{eqnarray*}
\int_X h(-\f_j) \om_{\f_j}^n &=& h(j) \int_{(\f \leq -j)} \om_{\f_j}^n
+\int_{(\f>-j)} h(-\f) \om_{\f_j}^n \\
&=& \int_{(\f \leq -j)} h(j) \om_{\f}^n+\int_{(\f>-j)} h(-\f) \om_{\f}^n \\
&\leq & \int_X h(-\f) \om_{\f}^n.
\end{eqnarray*}

Moreover if $B \subset X$ is a Borel subset, then
\begin{eqnarray*}
\left| \int_B h(-\f_j) \om_{\f_j}^n-\int_B h(-\f) \om_{\f}^n \right|
&\leq & 
\int_{B \cap (\f \leq -j)} h(-\f_j) \om_{\f_j}^n
+\int_{B \cap (\f \leq -j)} h(-\f) \om_{\f}^n \\
&\leq& 2 \int_{(\f \leq -j)} h(-\f) \om_{\f}^n \rightarrow 0.
\end{eqnarray*}
\end{proof}

\subsection{The comparison principle}

We now establish the {\it comparison principle}
which will be an important tool in the sequel.

\begin{thm}
Let $\f,\p \in {\mathcal E}(X,\om)$, then
$$
\int_{(\f<\p)} (\om+dd^c {\p})^n \leq \int_{(\f<\p)} (\om+dd^c {\f})^n.
$$

The class ${\mathcal E}(X,\om)$ is the largest
subclass of $PSH(X,\om)$ on which the operator 
$(\om+dd^c \cdot)^n$ is well defined and the comparison principle
is valid.
\end{thm}

\begin{proof}
Assume first that $\f,\p$ are bounded.
It follows from (1) that
\begin{eqnarray*}
\int_{(\f<\p)} \om_{\p}^n &=&
\int_{(\f<\p)} [\om+\max(\f,\p)]^n 
=\int_X \om^n-\int_{(\f \geq \p)} [\om+\max(\f,\p)]^n  \\
&\leq&  \int_X \om^n-\int_{(\f > \p)} \om_{\f}^n 
= \int_{(\f \leq \p)} \om_{\f}^n.
\end{eqnarray*}
Replacing $\p$ by $\p-\e$, $\e>0$, yields when $\e \searrow 0$,
$$
\int_{(\f<\p)} \om_{\p}^n=
\lim \nearrow  \int_{(\f <\p-\e)} \om_{\p}^n \leq 
\lim \nearrow  \int_{(\f  \leq \p-\e)} \om_{\f}^n
=\int_{(\f<\p)} \om_{\f}^n.
$$

When $\f,\p$ are unbounded, we set
$\f_j=\max(\f,-j)$ and $\p_k=\max(\p,-k)$. 
The comparison principle for bounded $\om$-psh functions
yields
$$
\int_{(\f_j < \p_k)} \om_{\p_k}^n \leq \int_{(\f_j<\p_k)} \om_{\f_j}^n.
$$
Observe that $(\f_j <\p) \subset (\f_j <\p_k) \subset (\f <\p_k)$.
Letting $k \rightarrow +\infty$ in the corresponding inequality yields,
by using monotone convergence theorem together with (3),
$$
\int_{(\f_j<\p)} \om_{\p}^n \leq \int_{(\f \leq \p)} \om_{\f_j}^n.
$$
Letting now $j \rightarrow +\infty$, we infer
$\int_{(\f<\p)} \om_{\p}^n \leq \int_{(\f \leq \p)} \om_{\f}^n$.
The desired inequality now follows by replacing $\p$
by $\p-\e$, $\e>0$, and letting
$\e \rightarrow 0^+$.

It turns out that the class ${\mathcal E}(X,\om)$ is the largest
class of $\om$-plurisubharmonic functions on which the complex
Monge-Amp\`ere operator is well defined and the comparison principle
is valid. Indeed let ${\mathcal F}$ be the largest
class with these properties, so that
${\mathcal E}(X,\om) \subset {\mathcal F} \subset PSH(X,\om)$
and
$$
\int_{(u \leq \p)} (\om+dd^c \p)^n \leq \int_{(u \leq \p)} (\om+dd^c u)^n,
$$
for all $u,\p \in {\mathcal F}$. Note that this inequality
is equivalent to the comparison
principle 1.5 since the measures $(\om+dd^c u)^n,(\om+dd^c \p)^n$ have
the same total mass $\int_X \om^n$. Fix $\f \in {\mathcal E}(X,\om)$
and $\p \in {\mathcal F}$ and apply previous inequality with
$u=\f+c$, $c \in \R$, to obtain
$$
(\om+dd^c \p)^n(\f=-\infty)=\lim_{c \rightarrow +\infty} 
\int_{(\f+c \leq \p)} (\om+dd^c \p)^n \leq \om_\f^n(\f=-\infty)=0.
$$
Since ${\mathcal E}(X,\om)$
characterizes pluripolar sets (Example 2.14), we infer
$\om_{\p}^n(\p=-\infty)=0$, hence $\p \in {\mathcal E}(X,\om)$,
i.e. ${\mathcal F}={\mathcal E}(X,\om)$.
\end{proof}

This principle allows us to derive 
important properties of the class
${\mathcal E}(X,\om)$.

\begin{pro}
The class ${\mathcal E}(X,\om)$ is convex.
Moreover if $\f \in {\mathcal E}(X,\om)$ and $\p \in PSH(X,\om)$
are such that $\f \leq \p$, then $\p \in {\mathcal E}(X,\om)$.
\end{pro}

\begin{proof}
The proof follows from the comparison principle
together with the following elementary observation:
{\it If $\f \in PSH(X,\om)$ and 
$\f/2 \in {\mathcal E}(X,\om)$, then $\f \in {\mathcal E}(X,\om)$}.
Indeed set $u=\f/2$, $u_j:=\max(u,-j)$, and $\f_j:=\max(\f,-j)$.
Observe that $u_j=\f_{2j}/2$ and
$$
\om_{u_j}:=\om+dd^c(\f_{2j}/2)=\frac{1}{2}(\om+\om_{\f_{2j}}) \geq 
\frac{1}{2}\om_{\f_{2j}},
$$
therefore
$$
\int_{(\f \leq -2j)} (\om+dd^c \f_{2j})^n
=\int_{(u \leq -j)} (\om+dd^c \f_{2j})^n 
\leq 2^n \int_{(u \leq -j)} (\om+dd^c u_j)^n \rightarrow 0.
$$

We now use the comparison principle to show 
that if $\f \in {\mathcal E}(X,\om)$, $\p \in PSH(X,\om)$
are such that $\f \leq \p$, then $\p/2 \in {\mathcal E}(X,\om)$,
so that $\p \in {\mathcal E}(X,\om)$ by previous observation.
Set $v:=\p/2$ and $v_j:=\max(v,-j)$, $\f_j:=\max(\f,-j)$.
We can assume without loss of generality that $\f \leq \p \leq -2$, 
hence $v \leq -1$.
It follows from $\f \leq \p$ that
$$
(v \leq -j) \subset (\f_{2j}<v_j-j+1) \subset (\f \leq -j),
$$
where the last inclusion simply uses $v_j \leq -1$.
We infer
$$
\om_{v_j}^n(v \leq -j) \leq \int_{(\f_{2j}<v_j-j+1)} \om_{v_j}^n
\leq \int_{(\f_{2j}<v_j-j+1)} \om_{\f_{2j}}^n
\leq \om_{\f_{2j}}^n(\f \leq -j) \rightarrow 0,
$$
as follows from Lemma 1.2. Thus $v=\p/2 \in {\mathcal E}(X,\om)$.

Using the comparison principle again, we now show that
if $\f,\p \in {\mathcal E}(X,\om)$, then $(\f+\p)/4 \in {\mathcal E}(X,\om)$.
It follows then from our first observation that $(\f+\p)/2 \in {\mathcal E}(X,\om)$,
thus ${\mathcal E}(X,\om)$ is convex.
Set $w:=(\f+\p)/4$, $w_j:=\max(w,-j)$, $\f_j:=\max(\f,-j)$ and $\p_j:=\max(\p,-j)$.
Observe that $(v \leq -j) \subset (\f \leq -2j) \cup (\p \leq -2j)$,
thus it suffices to show that 
$\om_{v_j}^n(\f \leq -2j) \rightarrow 0$. Assuming as above that
$\f,\p \leq -2$ yields
$(\f \leq -2j) \subset (\f_{2j}<w_j-j+1) \subset (\f \leq -j)$, hence
$$
\om_{v_j}^n(\f \leq -2j) \leq \om_{\f_{2j}}^n(\f \leq -j)
=\om_{\f_j}^n(\f \leq -j) \longrightarrow 0,
$$
as follows from Lemma 1.2.
\end{proof}

\begin{cor}
If $\f \in {\mathcal E}(X,\om)$, $\p \in PSH(X,\om)$, 
then $\max(\f,\p) \in {\mathcal E}(X,\om)$ and
$$
\1_{\{\f>\p\}} [\om+dd^c \f]^n=
\1_{\{\f>\p\}} [\om+dd^c \max(\f,\p)]^n.
$$
\end{cor}

\begin{proof}
It follows from the previous proposition that
$u:=\max(\f,\p)$ belongs to
${\mathcal E}(X,\om)$.
Set $\f_j=\max(\f,-j)$, $\p_{j+1}=\max(\p,-j-1)$
and $u_j:=\max(u,-j)$.
Observe that $\max(\f_j,\p_{j+1})=\max(\f,\p,-j)=u_j$.
Applying (1) yields
$$
\1_{\{\f_j>\p_{j+1}\}} (\om+dd^c \f_j)^n=
\1_{\{\f_j>\p_{j+1}\}} (\om+dd^c u_j)^n.
$$
Recall from (3) that 
$\1_{(\f>\p)} \om_{\f_j}^n \rightarrow \1_{(\f>\p)} \om_{\f}^n$.
Now $(\f>\p) \subset (\f_j>\p_{j+1})$ and 
$(\f_j>\p_{j+1}) \setminus (\f>\p) \subset (\f \leq -j)$, hence
$$
0 \leq \left[ \1_{(\f_j>\p_{j+1})}-\1_{(\f>\p)} \right] \om_{\f_j}^n
\leq \1_{(\f \leq -j)} \om_{\f_j}^n \rightarrow 0,
$$
since $\f \in {\mathcal E}(X,\om)$.
This shows that 
$$
\1_{\{\f_j>\p_{j+1}\}} (\om+dd^c \f_j)^n
\longrightarrow \1_{\{\f>\p\}} (\om+dd^c \f)^n.
$$
One proves similarly that
$\1_{\{\f_j>\p_{j+1}\}} (\om+dd^c u_j)^n \rightarrow 
\1_{\{\f>\p\}} (\om+dd^c u)^n$, observing
that $(\f_j>\p_{j+1}) \setminus (\f>\p) \subset (\max[\f,\p] \leq -j)$.
\end{proof}

The previous proposition also shows that
functions which belong to the class ${\mathcal E}(X,\om)$,
although possibly unbounded,
have mild singularities.

\begin{cor}
A function $\f \in {\mathcal E}(X,\om)$ 
has zero Lelong number at every point
$x \in X$.
\end{cor}

Note that this is not a sufficient condition to belong
to ${\mathcal E}(X,\om)$. It is for instance well known that, when $n=1$,
there are subharmonic functions whose Laplacian
has no Dirac mass but nevertheless charges a polar set.

\begin{proof}
For any point $x \in X$, one can construct a function $\r \in PSH(X,\om)$
which is smooth but at point $x$, and such that
$\r \sim c \log \text{dist}(\cdot,x)$ near $x$, for some constant $c>0$.
Such a function $\r$ has well defined Monge-Amp\`ere measure
$(\om+dd^c \f)^n$ and the latter has mass $\geq c^n$ at point $x$
(see [Dem]). Therefore $\rho \notin {\mathcal E}(X,\om)$.

Now if $\f \in PSH(X,\om)$ has a positive Lelong number at point $x$,
then $\f \leq \g \r+C$ on $X$, for some constants $\g,C>0$,
so it follows from the previous result that $\f$ does not belong to 
${\mathcal E}(X,\om)$.
\end{proof}

We will see in the next section (Theorem 2.6) that
the complex Monge-Amp\`ere operator
$\f \in {\mathcal E}(X,\om) \mapsto (\om+dd^c \f)^n$
is continuous under {\it any} decreasing sequences.

\begin{thm}
Let $\f_j \in {\mathcal E}(X,\om)$ be any sequence decreasing towards 
$\f \in {\mathcal E}(X,\om)$.
Then $(\om+dd^c {\f_j})^n \longrightarrow (\om+dd^c {\f})^n$.
\end{thm}

The proof of this result uses some properties of the weighted Monge-Amp\`ere 
energy that we introduce in section 2.
The following consequence will be quite useful when solving Monge-Amp\`ere
equations in section 4. It is due to E.Bedford and A.Taylor
[BT 1,3] when $\f,\p$ are bounded.

\begin{cor}
Assume $\f,\p \in {\mathcal E}(X,\om)$ are such that
$(\om+dd^c \f)^n \geq \mu$ and $(\om+dd^c \p)^n \geq \mu$ for some
positive Borel measure $\mu$ on $X$.
Then 
$$
[\om+dd^c \max(\f,\p)]^n \geq \mu.
$$
\end{cor}

\begin{proof}
Observe that $\max(\f,\p) \in {\mathcal E}(X,\om)$ by Proposition 1.6.
It follows from Corollary 1.7 that
$$
[\om+dd^c \max(\f,\p)]^n \geq \1_{\{\f>\p\}} (\om+dd^c \f)^n
+\1_{\{\f<\p\}} (\om+dd^c \p)^n 
\geq \1_{\{\f \neq \p\}} \mu.
$$
Thus we are done if $\mu(\f=\p)=0$.

Assume now that $\mu(\f=\p)>0$.
We show hereafter that $\mu(\f=\p+t_0)=0$ for
all $t_0 \in \R \setminus I_{\mu}$ where $I_{\mu}$ is at most countable.
Assuming this,
we can find a decreasing sequence $\e_j \searrow 0$ such that
$\mu(\f=\p+\e_j)=0$. Replacing $\p$ by $\p+\e_j$ above
yields
$$
[\om+dd^c \max(\f,\p+\e_j)]^n \geq \mu.
$$
The desired inequality therefore follows from Theorem 1.9.

It remains to show that $I_{\mu}:=\{ t_0 \in \R \, / \, \mu(\f=\p+t_0)>0 \}$
is at most countable. Consider
$
f:t \in \R \mapsto \mu(\{ \f<\p+t \}) \in \R^+.
$
This is an increasing function which is left continuous since $\mu$ is a Borel
measure. Moreover
$$
\lim_{t \stackrel{t>t_0}{\rightarrow} t_0} f(t)
=\mu\left( \{ \f \leq \p+t_0 \} \setminus \{ \p=-\infty\} \right)
=\mu(\{ \f \leq \p+t_0 \} )
$$
since $\mu(\p=-\infty) \leq \om_{\p}^n(\p=-\infty)=0$.
Therefore $f$ is continuous at $t_0$ unless
$\mu(\f=\p+t_0)>0$. Thus the set $I_{\mu}$ coincides with the set
of discontinuity of $f$ which is at most countable.
\end{proof}

\section{Low-energy classes}

Let $\chi:\R^- \rightarrow \R^-$ be a convex increasing function
such that $\chi(-\infty)=-\infty$.
It follows from the convexity assumption that
$$
0 \leq (-t) \chi'(t) \leq (-\chi)(t)
+\chi(0), \; \; \text{ for all } t \in \R^-.
$$

A straightforward computation shows that
$\chi \circ \f \in PSH(X,\om)$ whenever $\f \in PSH(X,\om)$ is non positive
and such that $\chi' \circ \f \leq 1$. Indeed
$$
dd^c \chi \circ \f=\chi'' \circ \f \, d\f \wedge d^c \f+\chi' \circ \f \, dd^c \f \geq -\om.
$$
Our aim here is to study the class of
$\om$-psh functions with finite $\chi$-energy.

\begin{defi}
We let
${\mathcal E}_{\chi}(X,\om)$ denote the set of 
$\om$-plurisubharmonic 
functions with
finite $\chi$-energy, i.e. 
$$
{\mathcal E}_{\chi}(X,\om):=\left\{ \f \in {\mathcal E}(X,\om) \, / \, 
 \chi(-|\f|) \in L^1( (\om+dd^c \f)^n) \right\}.
$$
\end{defi}

Observe that this definition is invariant under translation
both of the function $\f \mapsto \f+c$, $c \in \R$, and of the weight
$\chi \mapsto \chi-\chi(0)$.
We shall therefore often assume that the functions we are dealing
with are non-negative, and we will always normalize the weight
by requiring $\chi(0)=0$. We let
$$
{\mathcal W}^-:=\left\{ \chi:\R^- \rightarrow \R^- \, / \, 
\chi \text{ convex increasing, }
\chi(0)=0,\chi(-\infty)=-\infty
 \right\}
$$
denote the set of admissible weights.

\begin{pro}
$$
{\mathcal E}(X,\om)=\bigcup_{\chi \in {\mathcal W}^-} \E.
$$
\end{pro}

\begin{proof}
Fix $\f \in {\mathcal E}(X,\om)$.
Its complex Monge-Amp\`ere measure $(\om+dd^c \f)^n$ is well defined
and does not charge the pluripolar set $(\f=-\infty)$.
One can construct a continuous increasing function
$h:\R^+ \rightarrow \R^+$ such that
$h(+\infty)=+\infty$ and $h \circ |\f| \in L^1((\om+dd^c \f)^n)$.
Note that we can assume without loss of generality that
$h$ is concave (replacing if necessary 
$h$ by a concave increasing minorant $\tilde{h} \leq h$
such that $\tilde{h}(+\infty)=+\infty$).
Therefore $\f \in \E$, where
$\chi(t):=-h(-t)$.
\end{proof}

Observe that the union is increasing in the sense that
$\E \subset {\mathcal E}_{\tilde{\chi}}(X,\om)$
whenever $\tilde{\chi}=O(\chi)$ at infinity.
Observe also that $\chi=O(Id)$, as follows from the convexity assumption,
hence for all $\chi \in {\mathcal W}^-$,
$$
\E \supset {\mathcal E}^1(X,\om):=\{ \f \in {\mathcal E}(X,\om) \, / \, 
\f \in L^1(\om_{\f}^n) \}.
$$
When $n=\dim_{\C} X=1$, ${\mathcal E}^1(X,\om)$ is the classical class
of quasi-subharmonic functions of finite (unweighted) energy. 
\vskip.2cm

In the rest of section 2, $\chi \in {\mathcal W}^-$
denotes a fixed admissible weight and 
$$
E_{\chi}(\f):=\int_X (-\chi) \circ \f \, (\om+dd^c \f)^n
$$
denotes the $\chi$-energy of a function $\f \in \E$, $\f \leq 0$.

\subsection{Useful inequalities}

Many properties of the class $\E$ follow from simple integration by parts,
as shown by the following result.

\begin{lem}[{\bf The fundamental inequality}]
Let $\f,\p \in PSH(X,\om) \cap L^{\infty}(X)$ be such that $\f \leq \p \leq 0$.
Then
$$
0 \leq \int_X (-\chi) \circ \p \, (\om+dd^c \p)^n \leq 
2^n \int_X (-\chi) \circ \f \, (\om+dd^c \f)^n .
$$
\end{lem}

\begin{proof}
The proof follows from a repeated application of the following inequality, of independent
interest: let $T$ be any positive closed current of bidimension (1,1) on $X$, then
\begin{equation}
0 \leq \int_X (-\chi) \circ \f \, \om_{\p} \wedge T \leq 2 
\int_X (-\chi) \circ \f \, \om_{\f} \wedge T.
\end{equation}
Indeed observe that 
$\int_X (-\chi) \circ \p \, \om_{\p}^n \leq \int_X (-\chi) \circ \f \, \om_{\p}^n$,
hence it suffices to apply (6) with
$T=\om_{\f}^j \wedge \om_{\p}^{n-1-j}$, $0 \leq j \leq n-1$, to conclude.

It remains to prove (6). 
Here -- and quite often in the sequel -- we are going to use Stokes theorem
which yields
$\int_X udd^c v \wedge T=\int v dd^c u \wedge T$, whenever $u,v$ are bounded
$\om$-psh functions. Let us stress that there is no need for a {\it global}
regularization of $\om$-psh functions to justify this integration by parts.
We simply use the fact that $udd^c v \wedge T$ and $vdd^c u \wedge T$
are well defined and cohomologous currents on $X$, since
$$
udd^c v \wedge T-v dd^c u \wedge T=d [u d^c v \wedge T-v d^c u \wedge T].
$$
The latter computation can be justified by using local regularizations
together with continuity results of [BT 3]. We can also assume the weight
$\chi$ is smooth and then approximate it by using convolutions.
Observe that $\chi' \circ \f \, \om+dd^c (\chi \circ \f) \geq 0$, hence
\begin{eqnarray*}
\lefteqn{ \! \! \! \! \! \! \! \!
\! \! \! \! \! \!  
0 \leq \int_X (-\chi) \circ \f \, \om_{\p} \wedge T
=\int_X (-\chi) \circ \f \, \om \wedge T+
\int_X  (-\p) \, dd^c  \chi \circ \f \wedge T} \\
&\leq& \int_X (-\chi) \circ \f \, \om \wedge T+
\int_X  (-\f) \, \left[ \chi' \circ \f \, \om+dd^c  (\chi \circ \f )\right]\wedge T \\
&=& \int_X  (-\f) \chi' \circ \f \, \om \wedge T
+\int_X (-\chi) \circ \f \, \om_{\f} \wedge T.
\end{eqnarray*}
Observe now that $0 \leq (-\f) \chi' \circ \f \leq (-\chi) \circ \f$ and 
\begin{eqnarray*}
\int_X (-\chi) \circ \f \, \om_{\f} \wedge T &=&
\int_X (-\chi) \circ \f \, \om \wedge T
+\int_X \chi' \circ \f \, d \f \wedge d^c \f \wedge T \\
&\geq& \int_X (-\chi) \circ \f \, \om \wedge T,
\end{eqnarray*}
therefore
$\int_X (-\f) \chi' \circ \f \, \om \wedge T 
\leq \int_X (-\chi) \circ \f \, \om_{\f} \wedge T$,
which yields (6).
\end{proof}

It follows from Lemma 2.3 that the class $\E$ is stable
under taking maximum.
Another consequence of Lemma 2.3 is that
a given $\om$-psh function belongs to $\E$ if and only if {\it any} sequence
$\f_j \in PSH(X,\om) \cap L^{\infty}(X)$ decreasing to $\f$
has uniformly bounded $\chi$-energy. 

\begin{cor}
Fix $\f \in {\mathcal E}(X,\om)$.
The following conditions are equivalent:
\begin{enumerate}
\item $\f \in \E$;
\item for any sequence $\f_j \in PSH(X,\om) \cap L^{\infty}(X)$ decreasing towards $\f$,
$\sup_{j \geq 0} \int_X (-\chi)(|\f_j|) (\om+dd^c \f_j)^n <+\infty$;
\item there exists one sequence as in (2).
\end{enumerate}
\end{cor}

\begin{proof}
It follows from Proposition 1.4 that $\f \in \E$ if and only if
$\sup_{j \geq 0} \int_X (-\chi)(-|\f_j|) (\om+dd^c \f_j)^n <+\infty$,
where $\f_j:=\max(\f,-j)$ is the canonical approximating sequence.
Thus (1) is equivalent to (3).

Assume (3) holds. If $\p_j \in PSH(X,\om) \cap L^{\infty}(X)$
is any other sequence decreasing towards $\f$, 
then $\p_j \geq \f_{k_j}$ for some
(possibly large) $k_j \in \N$, hence by Lemma 2.3,
the sequence $\int_X (-\chi)(-|\p_j|) (\om+dd^c \p_j)^n$
is still uniformly bounded, showing (2).
The reverse implication $(2) \Rightarrow (3)$ is obvious.
\end{proof}

We now establish an inequality similar to (6), without any assumption on the 
relative localization of $\f,\p$.

\begin{pro}
Let $T$ be a positive closed current of bidimension (j,j) on $X$, $0 \leq j \leq n$,
and let $\f,\p \leq 0$ be bounded $\om$-psh functions. Then
$$
0 \leq \int_X (-\chi) \circ \f \, \om_{\p}^j \wedge T
\leq 2 \int_X (-\chi) \circ \f \, \om_{\f}^j \wedge T+
2\int_X (-\chi) \circ \p \, \om_{\p}^j \wedge T.
$$
\end{pro}

Observe that when $T$ has bidimension $(0,0)$, $T \equiv 1$,
this shows that 
$$
\chi \circ \E \subset L^1(\mu),
$$
for any Monge-Amp\`ere measure $\mu=(\om+dd^c \p)^n$, $\p \in \E$.
We shall use this fact in section 4, when describing the range of the 
complex Monge-Amp\`ere 
operator on classes $\E$.

\begin{proof}
Observe that $\chi'(2t) \leq \chi'(t)$ for all $t<0$, hence
$$
\int_X (-\chi) \circ \f \, \om_{\p}^j \wedge T=
\int_{-\infty}^0 \chi'(t) \om_{\p}^j \wedge T (\f <t) dt
\leq 2 \int_{-\infty}^0 \chi'(t) \om_{\p}^j \wedge T (\f <2t) dt.
$$
Now $(\f <2t) \subset (\f<\p+t) \cup (\p<t)$, hence
$$
\int_X (-\chi) \circ \f \, \om_{\p}^j \wedge T
\leq 2 \int_{-\infty}^0 \chi'(t) \om_{\p}^j \wedge T (\f <\p+t) dt
+2 \int_X (-\chi) \circ \p \, \om_{\p}^j \wedge T.
$$
The comparison principle yields
$\om_{\p}^j \wedge T (\f<\p+t) \leq \om_{\f}^j \wedge T(\f<\p+t)$.
The desired inequality follows by observing that $(\f<\p+t) \subset (\f<t)$.
\end{proof}

\subsection{Continuity of weighted Monge-Amp\`ere operators}

We are now in position to prove a strong version of Theorem 1.9.

\begin{thm}
Let $\f_j \in PSH(X,\om)$ be a sequence decreasing towards $\f \in \E$.
Then $\f_j \in \E$ and
$(\om+dd^c {\f_j})^n \longrightarrow (\om+dd^c{\f})^n$.
Moreover for any $\tilde{\chi} \in {\mathcal W}^-$
such that $\tilde{\chi}=o(\chi)$, one has
$$
\tilde{\chi} (-|\f_j|) \, (\om+dd^c \f_j)^n 
\longrightarrow \tilde{\chi} (-|\f|) \, (\om+dd^c \f)^n.
$$

When $\f_j:=\max(\f,-j)$ is the canonical approximation, then 
$$
\1_B \chi (-|\f_j|) (\om+dd^c \f_j)^n \rightarrow \1_B \chi (-|\f|) (\om+dd^c \f)^n
$$
for all Borel subsets $B \subset X$.
\end{thm}

Observe that Theorem 1.9 easily follows  from this result
together with  Proposition 2.2.

\begin{proof}
\noindent {\bf Step 1.}
Assume first that $\f_j:=\max(\f,-j)$.
For $B \subset X$, 
\begin{eqnarray*}
\left| \int_B |\chi| \circ \f_j \, \om_{\f_j}^n
-\int_B |\chi| \circ \f \, \om_{\f}^n \right| 
&\leq& \int_{B \cap (\f\leq -j)} 
\left[ |\chi| \circ \f_j \, \om_{\f_j}^n
+|\chi| \circ \f \, \om_{\f}^n \right] \\
& \leq & 2 \int_{(\f \leq -j)} |\chi| \circ \f \, \om_{\f}^n 
\rightarrow 0,
\end{eqnarray*}
since $\chi \circ \f \in L^1(\om_{\f}^n)$.
We have used here the upper bound
\begin{eqnarray*}
\int_{(\f \leq -j)} (-\chi) \circ \f_j \, \om_{\f_j}^n
&=&(-\chi)(-j) \int_{(\f \leq-j)} \om_{\f_j}^n \\
&=&(-\chi)(-j) \int_{(\f \leq-j)} \om_{\f}^n
\leq \int_{(\f \leq -j)} (-\chi) \circ \f \, \om_{\f}^n.
\end{eqnarray*}
This shows that 
$\1_B \chi \circ \f_{\j} \, \om_{\f_j}^n \rightarrow \1_B \chi \circ \f \, \om_{\f}^n$.
We infer that the fundamental inequality holds when $\f,\p \in \E$.
\vskip.2cm

\noindent {\bf Step 2.} 
We now consider the case of a general sequence $(\f_j)$ that decreases
towards $\f$.
The continuity  of the complex Monge-Amp\`ere operator $(\om+dd^c \cdot)^n$
along such sequences is due to E.Bedford and A.Taylor [BT 3] when 
$\f$ is bounded, and we shall reduce the problem to this case.
We can assume without loss of generality that $\f,\f_j \leq 0$.
Consider
$$
\f_j^K:=\max(\f_j,-K) \; \text{ and } \f^K:=\max(\f,-K).
$$
The integer $K$ being fixed, the sequence $(\f_j^K)_j$ is uniformly
bounded and decreases towards $\f^K$, hence
$$
\left(\om+dd^c \f_j^K \right)^n \stackrel{j \rightarrow +\infty}{\longrightarrow} 
\left(\om+dd^c \f^K\right)^n.
$$
Thus we will be done if we can show that
$(\om+dd^c \f_j^K)^n$ converges towards $(\om+dd^c \f_j)^n$
as $K \rightarrow +\infty$, uniformly  with respect to $j$.
Let $h$ be a continuous test function on $X$. Then
\begin{eqnarray*}
\lefteqn{
\left| \langle h, (\om+dd^c \f_j^K)^n-(\om+dd^c \f_j)^n \rangle \right|} \\
& \leq&  ||h||_{L^{\infty}(X)}
\int_{(\f_j \leq -K)} \left[ (\om+dd^c \f_j^K)^n+(\om+dd^c \f_j)^n \right] \\
&\leq& \frac{||h||_{L^{\infty}(X)}}{(-\chi)(-K)}
\left\{ \int_X (-\chi) \circ \f_j^K (\om+dd^c \f_j^K)^n
+\int_X (-\chi) \circ \f_j (\om+dd^c \f_j)^n \right\}.
\end{eqnarray*}
Since $\f \leq \f_j \leq \f_j^K$, it follows now from Step 1 that
the last two integrals are uniformly bounded from above
by $2^n \int_X (-\chi) \circ \f \, \om_{\f}^n$.
This yields the desired uniformity.
Note that the same proof shows that
$\tilde{\chi} \circ \f_j \om_{\f_j}^n \, \rightarrow \tilde{\chi} \circ \f \, \om_{\f}^n$,
whenever  $\tilde{\chi}=o(\chi)$, so that a factor
$\tilde{\chi}(-K)/\chi(-K) \rightarrow 0$ yields uniformity.
\end{proof}

\begin{cor}
The fundamental inequality holds when $\f,\p \in \E$.
Moreover if $0 \geq \f_j \in \E$ is a sequence of functions 
converging towards $\f$ in $L^1(X)$ and such that
$E_{\chi}(\f_j)$ is uniformly bounded, then $\f \in \E$.
\end{cor}

\begin{proof}
Set $\Phi_j:=(\sup_{l \geq j} \f_l)^*$.
This is a sequence of $\om$-psh functions which decrease towards $\f$.
Since $0 \geq \Phi_j \geq \f_j$, it follows from Lemma 2.3 and Theorem 2.6
that $E_{\chi}(\Phi_j)$ is uniformly bounded. Thus $\f \in \E$
by Corollary 2.4.
\end{proof}

The class ${\mathcal E}^1(X,\om)$ plays a special role in many respects,
as will become clear in section 3.
We indicate here one specific
property that will be useful
when solving Monge-Amp\`ere equation in section 4.

\begin{pro}
Let $\f_j,\f \in {\mathcal E}^1(X,\om)$ be such that
$\f_j \rightarrow \f$ in $L^1(X)$.
 If $\int_X |\f_j-\f| \om_{\f_j}^n \rightarrow 0$, then
$(\om+dd^c {\f_j})^n \rightarrow (\om+dd^c{\f})^n$.
\end{pro}

\begin{proof}
We can assume without loss of generality 
that $\f_j,\f \leq 0$.
Passing to a subsequence if necessary, 
we can assume $\int |\f_j-\f| \om_{\f_j}^n \leq 1/j^2$.
Consider
$$
\Phi_j:=\max (\f_j,\f-1/j) \in {\mathcal E}^1(X,\om).
$$
It follows from Hartogs' lemma that $\Phi_j \rightarrow \f$ in capacity.
This means that
$Cap_{\om}(|\Phi_j-\f|>\e) \rightarrow 0$, for all $\e>0$,
where
$$
Cap_{\om}(K):=\sup \{ \int_K (\om+dd^c u)^n \, / \, u \in PSH(X,\om), \, -1 \leq u \leq 0 \}
$$
is the Monge-Amp\`ere capacity (see [BT 3], [GZ 1]).
It is a well-known consequence of the quasicontinuity of $\om$-psh functions
that $(\om+dd^c \Phi_j)^n \rightarrow (\om+dd^c \f)^n$,
when the $\Phi_j'$s are uniformly bounded [X]. We can reduce to this
case by showing that
$(\om+dd^c \Phi_j^K)^n$ converges towards $(\om+dd^c \f^K)^n$ uniformly with respect
to $K$, where $\Phi_j^K:=\max(\Phi_j,-K)$ and $\f^K:=\max(\f,-K)$.
Indeed if $\theta$ is a test function, then
\begin{eqnarray*}
\left| \langle \om_{\Phi_j^K}^n,\theta \rangle
- \langle \om_{\Phi_j}^n,\theta \rangle \right|
& \leq & \sup_X |\theta|
\left\{ \int_{(\Phi_j \leq -K)} \om_{\Phi_j^K}^n+
\int_{(\Phi_j \leq -K)} \om_{\Phi_j}^n \right\} \\
& \leq & \frac{\sup_X |\theta|}{K} 
\left\{ E_1(\Phi_j^K)+E_1(\Phi_j) \right\} \\
& \leq & 2 \frac{\sup_X |\theta|}{K} E_1 (\f-1),
\end{eqnarray*}
where $E_1=E_{\chi}$ for $\chi(t)=t$.
This follows from Lemma 2.3 and the lower bound $\Phi_j \geq \f-1$.
Thus we have shown that the measures $(\om+dd^c \Phi_j)^n$
converge towards $(\om+dd^c \f)^n$.

We now need to compare $(\om+dd^c {\Phi_j})^n$ and $(\om+dd^c {\f_j})^n$. 
It follows from Corollary 1.7, that
$$
\om_{\Phi_j}^n \geq \1_{\{ \f_j \geq \f-1/j \}} \cdot \om_{\f_j}^n.
$$
Let $E_j$ denote the set $X \setminus \{\f_j \geq \f-1/j\}$, i.e.
$E_j=\{\f-\f_j>1/j\}$. Our assumption implies that 
$\1_{E_j} \om_{\f_j}^n \rightarrow 0$, indeed
$$
0 \leq \int_{E_j} \om_{\f_j}^n \leq j \int_X |\f-\f_j| \om_{\f_j}^n \leq \frac{1}{j}.
$$
Therefore $0 \leq \om_{\f_j}^n \leq \om_{\Phi_j}^n+o(1)$, hence
$\om_{\f}^n=\lim \om_{\f_j}^n$.
\end{proof}

\begin{rqe}
Observe the following property which was used and proved above (when $\chi(t)=t$):
if $\f \in PSH(X,\om)$, $\f_j \in {\mathcal E}(X,\om)$,
$\f_j \rightarrow \f$ in capacity,
and $\f_j \geq \p$ for some fixed $\p \in {\mathcal E}(X,\om)$,
then $\f \in {\mathcal E}(X,\om)$ and 
$(\om+dd^c \f_j)^n \rightarrow (\om+dd^c \f)^n$.
\end{rqe}

\subsection{Homogeneous weights}

The weights $\chi_p(t)=-(-t)^p$, $0<p \leq 1$,
belong to ${\mathcal W}^-$.
We shall use the notation
$$
{\mathcal E}^p(X,\om):={\mathcal E}_{\chi}(X,\om),
\text{ when } \chi(t)=-(-t)^p.
$$
These classes are easier to understand thanks to the homogeneity
property of the weight function, $\chi(\e t)=\e^p \chi(t)$.

Also the class ${\mathcal E}^1(X,\om)$ deserves special attention, as
it is the turning point between low-energy and high-energy classes.
All functions from ${\mathcal E}^1(X,\om)$ have gradient in $L^2(X)$
(see Proposition 3.2),
while most functions of lower energy do not (Example 2.14).
We will show (Theorem 3.3) that solutions of complex Monge-Amp\`ere
equations have a unique solution in ${\mathcal E}^1(X,\om)$,
while it is a question that remains open in classes of lower energy.
\vskip.2cm

Our aim here is to establish further properties of the
classes ${\mathcal E}^p(X,\om)$. It will allow us to
give a complete characterization of the range of the
complex Monge-Amp\`ere operator on them (see Theorem 4.2).

\begin{pro}
There exists $C_p>0$ such that for all 
$0 \geq  \f_0,\ldots,\f_n \in PSH(X,\om) \cap L^{\infty}(X)$,
$$
0 \leq \int_X (-\f_0)^p \om_{\f_1} \wedge \cdots \wedge \om_{\f_n}
\leq C_p \max_{0 \leq j \leq n} \left[ \int_X (-\f_j)^p \om_{\f_j}^n \right].
$$

In particular the class ${\mathcal E}^p(X,\om)$ is starshaped and convex.
\end{pro}

\begin{proof}
It follows from Proposition 2.5 applied with
$\f=\f_0$, $\p=\f_1$ and $T=\om_{\f_2} \wedge \cdots \wedge \om_{\f_n}$ that
\begin{equation}
\int_X (-\chi) \circ \f_0 \, \om_{\f_1} \wedge T
\leq 2 \int_X (-\chi) \circ \f_0 \, \om_{\f_0} \wedge T
+ 2 \int_X (-\chi) \circ \f_1 \, \om_{\f_1} \wedge T,
\end{equation}
thus we can assume $\f_0=\f_1$ in the sequel.

Set $u=\e \sum_{i=1}^n \f_i$, where $\e>0$ is small enough
($\e<(2n)^{-1/p}$ will do). Observe that
$\om_u^n \geq \e^n \om_{\f_1} \wedge \cdots \wedge \om_{\f_n}$, hence it suffices
to get control on $\int_X (-\chi) \circ \f_i \, \om_u^n$
for all $1 \leq i \leq n$ to conclude. 
By using Proposition 2.5 again,
$$
\int_X (-\chi) \circ \f_i \, \om_u^n \leq
2 E_{\chi}(\f_i)+2E_{\chi}(u),
$$
where $E_{\chi}(u):=\int_X (-\chi) \circ u \, \om_u^n$.
By subadditivity and homogeneity of $-\chi$, we get
$E_{\chi}(u) \leq \e^p \sum_{j=1}^n \int_X (-\chi) \circ \f_j \, \om_u^n$,
hence 
\begin{equation}
\sum_{i=1}^n \int_X (-\chi) \circ \f_i \, \om_u^n
\leq \frac{2}{1-2n\e^p} \sum_{i=1}^n E_{\chi}(\f_i).
\end{equation}
We deduce from (7), (8) and 
$\om_u^n \geq \e^n \om_{\f_1} \wedge \cdots \wedge \om_{\f_n}$
that
$$
0 \leq \int_X (-\f_0)^p \om_{\f_1} \wedge \cdots \wedge \om_{\f_n}
\leq \frac{4n}{\e^n [1-2n\e^p]} \max_{1 \leq i \leq n} E_{\chi}(\f_i).
$$

Each class $\E$ is clearly starshaped. It follows from previous 
inequality that the classes ${\mathcal E}^p(X,\om)$
are also convex.
\end{proof}

It follows from this result that if some functions $\p_j \in {\mathcal E}^p(X,\om)$,
$\p_j \leq 0$, have uniformly bounded energy $\sup_{j \geq 0 }E_{\chi}(\p_j) <+\infty$,
then
$$
\p:=\sum_{j \geq 1} 2^{-j} \p_j \in {\mathcal E}^p(X,\om).
$$
This observation, together with the homogeneity property of $\chi(t)=-(-t)^p$
allows us to derive the following quantitative characterization of
integrability properties with respect to a given measure $\mu$.

\begin{lem}
Let $\mu$ be a probability measure on $X$.
Then ${\mathcal E}^p(X,\om) \subset L^p(\mu)$ if and only if there
exists $C>0$ such that
for for all $\p \in PSH(X,\om) \cap L^{\infty}(X)$
with $\sup_X \p=-1$,
\begin{equation}
0 \leq \int_X (-\p)^p \, d\mu \leq C
\left[ \int_X (-\p)^p \, \om_{\p}^n \right]^{\frac{p}{p+1}}.
\end{equation}
\end{lem}

\begin{proof}
Assume on the contrary that (9) is not satisfied, then there exists
$\p_j\in PSH(X,\om) \cap L^{\infty}(X)$, $\sup_X \p_j=-1$, such that
$$
\int_X (-\chi) \circ \p_j \, d\mu \geq 4^{jp}
M_j^{\frac{p}{p+1}},
\text{ where } 
M_j:=E_{\chi}(\p_j).
$$

If $(M_j)$ is bounded we consider $\p:=\sum_{j \geq 1} 2^{-j} \p_j$.
This is a $\om$-psh function which belongs to
${\mathcal E}^p(X,\om)$ by Proposition 2.10. Now
$$
\int_X (-\chi) \circ \p d\mu \geq \int_X (-\chi)(2^{-j} \p_j) d\mu
\geq 2^{jp} M_j^{\frac{p}{p+1}} \geq 2^{jp},
$$
because $\p_j \leq -1$, hence $M_j \geq 1$. Thus
$\int_X (-\chi) \circ \p d\mu=+\infty$, a contradiction.

We similarly obtain a contradiction if $(M_j)$ admits a bounded
subsequence, so we can assume $M_j \rightarrow +\infty$ and $M_j \geq 1$.
Set $\f_j=\e_j \p_j$, where $\e_j=M_j^{-1/(p+1)}$. We show herebelow that
$E_{\chi}(\f_j)$ is bounded, hence
$\f:=\sum_{j \geq 1} 2^{-j} \f_j \in {\mathcal E}^p(X,\om)$
by Proposition 2.10. Now
\begin{eqnarray*}
\int_X (-\chi) \circ \f \, d\mu 
& \geq & 2^{-jp} \int_X (-\chi) \circ \f_j \, d\mu 
\geq 2^{-jp}\e_j^p \int_X (-\chi) \circ \p_j \, d\mu \\
& \geq & 2^{jp} \e_j^p M_j^{p/(p+1)}=2^{jp},
\end{eqnarray*}
hence $\int_X (-\chi) \circ \f \, d\mu =+\infty$, a contradiction.

It remains to check that $(E_{\chi}(\f_j))$ is indeed bounded.
Observe that $\om_{\f_j} \leq \e_j \om_{\p_j}+\om$, thus
\begin{eqnarray*}
\lefteqn{ E_{\chi}(\f_j) =\e_j^p \int_X (-\chi) \circ \p_j \, \om_{\f_j}^n } \\
&& \leq \e_j^p \left[ \int_X (-\chi) \circ \p_j \, \om^n+2^n \e_j E_{\chi}(\p_j) \right] \\
&& \leq \int_X (-\p_j) \om^n+2^n.
\end{eqnarray*}
We have used here the definition of $\e_j$,
$\e_j^{p+1} E_{\chi}(\p_j)=1$, the bounds 
$\e_j \leq 1$,  $|\chi|(t) \leq |t|$, and the inequalities
\begin{eqnarray*}
E_{\chi}(\p_j)&=&\int_X (-\chi) \circ \p_j \, \om \wedge \om_{\p_j}^{n-1}+
\int_X \chi' \circ \p_j \, d\p_j \wedge d^c \p_j \wedge \om_{\p_j}^{n-1} \\
&\geq & \int_X (-\chi) \circ \p_j \, \om \wedge \om_{\p_j}^{n-1}
\geq \int_X (-\chi) \circ \p_j \, \om^k \wedge \om_{\p_j}^{n-k},
\end{eqnarray*}
for all $1 \leq k \leq n-1$.
Now $\int_X (-\p_j) \om^n$ is bounded because $\sup_X \p_j=-1$ (see Proposition 2.7 in [GZ 1]),
hence $(E_{\chi}(\f_j))$ is bounded.
\end{proof}

\begin{rqe}
These results actually apply for any weight
$\chi \in {\mathcal W}^-$ close enough to a homogeneous weight,
as we indicate in section 3.3.
\end{rqe}

\subsection{Some examples}

The function $L(t):=-\log(1-t)$ belongs to ${\mathcal W}^-$,
as well as $L_p(t):=L \circ \cdots \circ L(t)$ ($p$ times).
It is in fact necessary to consider functions $\chi$ with 
arbitrarily slow growth in order to understand the range
of the complex Monge-Amp\`ere operator on
${\mathcal E}(X,\om)$, as the following example shows.

\begin{exa}
Fix $h \in {\mathcal W}^- \cap {\mathcal C}^{\infty}(\R^-,\R^-)$
with $h'(-\infty)=0$, and consider
$$
\f: z \in \C \subset \P^1=\C \cup \{\infty \}
\mapsto h \left( \log|z|-\frac{1}{2}\log[1+|z|^2]-1\right)
\in \R^-,
$$
with the convention $\f(\infty)=h(-1)$.
This is an $\om$-subharmonic function on the Riemann sphere $\P^1$,
where $\om$ denotes the Fubini-Study volume form.
Note that $\f$ is smooth but at the origin $0 \in \C$ and 
$$
\om_{\f}=[f(z)+o(f(z))] dV, \; \; 
f(z)=\frac{c}{|z|^2}h''(\log|z|), \text{ near 0},
$$
where $dV$ denotes the euclidean volume form and $c>0$.
One thus checks that 
$\f \in \E$ if $\chi$ does not grow too fast
(e.g. $|\chi|(t) \lesssim [h' \circ h^{-1}(t)]^{-1/2}$),
while $\chi \circ \f \notin L^1(\om_{\f})$
(hence $\f \notin \E$) if e.g. $|\chi(t)| \geq [h' \circ h^{-1}(t)]^{-1}$.

As a concrete example, consider $h(t)=t\log(1-t)$.
The reader can check that $h$ satisfies our requirements, and that
in this case $\om_{\f}=f dV$ is a measure with density such that
$$
f(z) \sim \frac{c'}{|z|^2 (-\log |z|^2) [\log (-\log|z|^2)]^2}
\text{ near } 0.
$$
In particular $|\f|^p \notin L^1(\om_{\f})$ for any $p>0$, moreover
$\log[1-\f] \notin L^1(\om_{\f})$.
One can obtain explicit examples with even slower growth by 
considering $h(t)=L^p(t):=L \circ \cdots \circ L(t)$ ($p$ times),
where $L(t)=-\log(1-t)$.
\end{exa}

Our next example shows that it is possible to define and control the 
complex Monge-Amp\`ere measure of functions with 
slightly attenuated singularites, although these functions need not have
gradient in $L^2(X)$.

\begin{exa}
For $\f \in PSH(X,\om)$, 
$\f \leq -1$, and $0<p<1$, we set $\f_p:=-(-\f)^p$.
Then
$$
\om_{\f_p}=p(-\f)^{p-1} \om_{\f}+[1-p(-\f)^{p-1}]\om+p(1-p)(-\f)^{p-2}d\f \wedge d^c \f \geq 0,
$$
hence $\f_p \in PSH(X,\om)$. One can compute similarly
$(\om+dd^c \f_p)^n$ and check that there exists $\a(p,n)>0$ small enough
and $C(p,n)>0$ independent of $\f$ such that
$$
\int_X (-\f_p)^{\a(p,n)} (\om+dd^c \f_p)^n \leq C(p,n).
$$
This shows that $\f_p \in \E$ with $\chi(t)=-(-t)^{\a(p,n)}$, as soon as $p<1$.
In particular the complex Monge-Amp\`ere operator
of $\f_p=-(-\f)^p$ is well defined for any $\om$-psh function $\f \leq -1$ and
for any $p<1$, although $\nabla \f_p$ does not belong to $L^2(X)$
when $\om_{\f}$ is the current of integration along a complex hypersurface
and $1/2 \leq p <1$.
\end{exa}

\section{High energy classes}

We now consider further classes of $\om$-plurisubharmonic with
milder singularities. Not only are they interesting in themselves,
but we also actually need first to understand the range of
the complex Monge-Amp\`ere operator on them, before being able
to describe the corresponding range on previous classes (see
the proofs of Theorems 4.1, 4.2). 

We let
${\mathcal W}^+$ denote the set of 
{\it concave} increasing functions $\chi:\R^- \rightarrow \R^-$
such that  $\chi(0)=0,\chi(-\infty)=-\infty$ and we set
$$
{\mathcal W}^+_M:=\left\{ \chi \in {\mathcal W}^+ \, / \, 
0 \leq |t \chi'(t)| \leq M |\chi(t)|, \;
\text{ for all } t \in \R^- \right\},
$$
where $M>0$ is independent of $t$.
For $\chi \in {\mathcal W}^+$ we consider
$$
\E:=\{ \f \in {\mathcal E}(X,\om) \, / \, \chi (-|\f|) \in L^1((\om+dd^c \f)^n) \}.
$$

Observe that for all weights $\chi_1 \in {\mathcal W}^-, \chi_2 \in {\mathcal W}^+$,
we have
$$
PSH(X,\om) \cap L^{\infty}(X) \subset 
{\mathcal E}_{\chi_2}(X,\om) \subset {\mathcal E}^1(X,\om)
\subset {\mathcal E}_{\chi_1}(X,\om)
\subset {\mathcal E}(X,\om).
$$
Thus the class ${\mathcal E}^1(X,\om)$ is the turning point
between classes of low energy and those of high energy.

\subsection{Gradient of quasiplurisubharmonic functions}

The classes $\E$, $\chi \in {\mathcal W}:={\mathcal W}^- \cup {\mathcal W}^+$,
form a whole scale among unbounded $\om$-plurisubharmonic functions,
joining the maximal class ${\mathcal E}(X,\om)$
of $\om$-psh functions $\f$ whose Monge-Amp\`ere measure $(\om+dd^c \f)^n$ 
does not charge pluripolar sets (Proposition 2.2) to the class
$PSH(X,\om) \cap L^{\infty}(X)$ of bounded $\om$-psh functions.

\begin{pro}
$$
PSH(X,\om) \cap L^{\infty}(X)=\bigcap_{\chi \in {\mathcal W}^+} \E.
$$
\end{pro}

\begin{proof}
Clearly a bounded $\om$-psh function belongs to $\E$ for all 
$\chi \in {\mathcal W}$.

Conversely assume $\f \in \E$ for any weight $\chi \in {\mathcal W}$.
We claim there exists $t_{\infty}>0$ such that 
$(\om+dd^c \f)^n(\f<-t_{\infty})=0$.
Otherwise we could set
$$
\chi(t):=-h(-t),
\text{ where } h'(t)=\frac{1}{\om_{\f}^n(\f<-t)}
$$
is well defined for all $t>0$.
Now $\chi \in {\mathcal W}^+$ and, assuming $\f \leq 0$, 
$$
\int_X (-\chi) \circ \f \om_{\f}^n=\int_0 ^{+\infty} h'(t) \om_{\f}^n( \f<-t) dt=+\infty,
$$
contradicting $\f \in \E$.

We infer $(\om+dd^c \f)^n(\f < -t)=0$ for all $t \geq t_{\infty}$.
It follows from the comparison principle that
$Vol_{\om}(\f<-t_{\infty}) \leq \int_{(\f<-t_{\infty})} \om_{\f}^n$.
Thus $\f \geq -t_{\infty}$ almost everywhere,
hence everywhere.
\end{proof}

One of the main differences between functions with finite
low-energy and those with finite high-energy resides in
the following observation.

\begin{pro}
Fix $\chi \in {\mathcal W}^+$.
Then any function $\f \in {\mathcal E}_{\chi}(X,\om)$ is such that
$\nabla \f \in L^2(\om^n)$.
\end{pro}

\begin{proof}
It suffices to establish this result when $\chi(t)=t$,
i.e. when $\f \in {\mathcal E}^1(X,\om)$.
Let $T$ be a positive current of bidimension $(1,1)$ and fix
$\f \in PSH(X,\om) \cap L^{\infty}(X)$ such that $\f \leq 0$. Then
$$
\int_X (-\f) \om_{\f} \wedge T=\int_X (-\f) \om \wedge T +\int_X d\f \wedge d^c \f \wedge T
\geq \int_X (-\f) \om \wedge T.
$$
A repeated application of this inequality, applied with $T=\om_{\f}^j \wedge \om^{n-1-j}$,
therefore yields 
\begin{equation}
0 \leq \int_X d\f \wedge d^c \f \wedge \om^{n-1}
\leq \int_X (-\f) \om_{\f} \wedge \om^{n-1} \leq \int_X (-\f) \om_{\f}^{n}.
\end{equation}

Fix now $\f \in {\mathcal E}^1(X,\om)$ and set
$\f_j:=\max(\f,-j)$. We can assume without loss of generality
$\f,\f_j \leq 0$. It follows from (10) that the sequence
of gradients $(\nabla \f_j)$ has uniformly bounded $L^2$-norm.
Since $\nabla \f_j$ converges towards $\nabla \f$ in the weak sense
of distributions, it follows from the weak $L^2$-compactness
of $\nabla \f_j$ that $\nabla \f \in L^2(X)$.
\end{proof}

We can prove similarly that $\f \in {\mathcal E}_{\chi}(X,\om)$
satisfies all intermediate local boundedness conditions
required in order for the local complex Monge-Amp\`ere operator
$(\om+dd^c \f)^n$ to be well defined (see [Bl 3]).
It follows however from Example 2.14
that there are functions $\f \in \E$ such that $\nabla \f \notin L^2(X)$
when $\chi \in {\mathcal W}^-$.

\subsection{Extension of Calabi's uniqueness result}

When $\f,\p$ are smooth functions such that
$\om_{\f},\om_{\p}$ are K\"ahler forms, it has been shown by
E.Calabi [Ca] that 
$$
(\om+dd^c \f)^n=(\om+dd^c \p)^n \Longrightarrow \f-\p \equiv \text{constant.}
$$
Calabi's proof has been generalized by Z.Blocki [Bl 1] to the case
where $\f,\p \in PSH(X,\om) \cap L^{\infty}(X)$.
In both cases,  the proof consists
in showing that $\nabla(\f-\p)=0$,
by using ingenious integration by parts.

We push this argument further by showing that there is
still uniqueness in the class
${\mathcal E}^1(X,\om)$ (Theorem B in the introduction). 
We rely heavily on the fact that
$\nabla \f \in L^2(\om^n)$ when $\f \in {\mathcal E}^1(X,\om)$.
It is an interesting open question to establish the
uniqueness of solutions in classes of lower energy.

\begin{thm} 
Assume $(\om+dd^c{\f})^n \equiv (\om+dd^c{\p})^n$, where
$\f \in {\mathcal E}(X,\om)$ and $\p \in {\mathcal E}^1(X,\om)$.
Then $\f-\p$ is constant.
\end{thm}

\begin{proof}
We assume first that both $\f$ and $\p$ are in
${\mathcal E}^1(X,\om)$. 
Set $f=(\f-\p)/2$ and $h=(\f+\p)/2$. It follows from
Proposition 2.10 that $h \in {\mathcal E}^1(X,\om)$.
We can assume w.l.o.g. $\f,\p \leq -C_{\om}$, where 
$C_{\om}>0$ is chosen so that
$\int_X (-h)\om_h^n \geq 1$.
We are going to prove that $\nabla f=0$ by 
showing that $\int_X df \wedge d^c f \wedge \om^{n-1}=0$.
This will be done by establishing an upper bound 
involving the $E_1$-energy of $h$ and by using the
following formal computation,
$$
0 \leq \int_X df \wedge d^c f \wedge \om_h^{n-1}
\leq \int_X df \wedge d^c f \wedge \sum_{l=0}^{n-1} \om_{\f}^l \wedge \om_{\p}^{n-1-l}
=\int_X \frac{-f}{2} (\om_{\f}^n-\om_{\p}^n)=0.
$$

Of course we should (and could) justify these integration by parts and make sense
of all terms involved in this computation.
We are rather going to establish the following a priori bound,
when $\f,\p$ are bounded,
$$
{\bf (\dag)} \hskip.5cm
\int_X df \wedge d^c f \wedge \om^{n-1}
\leq 3^{n-1} \left( \int_X df \wedge d^c f \wedge \om_h^{n-1} \right)^{1/2^{n-1}}
\int_X (-h) \om_h^n.
$$
Approximating $\f,\p$ by $\f_j=\max(\f,-j), \p_j=\max(\p,-j)$, it will then 
follow from Theorem 2.6 and $(\dag)$ that
$2f=\f-\p$ is constant, if $\om_{\f}^n \equiv \om_{\p}^n$.

The a priori bound $(\dag)$ follows by
applying inductively the following inequality, where
$T=\om_h^l \wedge \om^{n-2-l}$ 
is a closed positive current of bidimension $(2,2)$
and $l=n-2,\ldots,0$,
$$
{\bf (\dag)}_T \hskip.1cm
\int_X df \wedge d^c f \wedge \om \wedge T \leq 
3 \left( \int_X df \wedge d^c f \wedge \om_h \wedge T \right)^{\frac{1}{2}}
\left( \int_X (-h) \om_h^2 \wedge T \right)^{\frac{1}{2}}.
$$
Indeed observe that for $T=\om_h^l \wedge \om^{n-2-l}$,
$$
\int_X (-h) \om_h^{2+l} \wedge \om^{n-2-l}
\leq \int_X (-h) \om_h^n,
\text{ for all } l=n-2,\ldots,0
$$
and use
that $\int_X (-h)\om_h^n \geq 1$ since $h \leq -C_{\om}$,
to derive $(\dag)$ from $(\dag)_T$.

We now establish $(\dag)_T$. Note that 
$$
df \wedge d^c f \wedge \om=df \wedge d^c f \wedge \om_h -
df \wedge d^c f \wedge dd^c h,
$$
hence integrating by parts yields
$$
\int_X df \wedge d^c f \wedge \om \wedge T=\int_X df \wedge d^c f \wedge \om_h \wedge T
+\int_X df \wedge d^c h \wedge \frac{(\om_{\f}-\om_{\p})}{2} \wedge T.
$$
It follows from the Cauchy-Schwarz inequality that
$$
\left| \int_X df \wedge d^c h \wedge \om_{\f} \wedge T \right|
\leq 2 \left(\int_X df \wedge d^c f \wedge \om_h \wedge T \right)^{\frac{1}{2}}
\cdot
\left(\int_X dh \wedge d^c h \wedge \om_h \wedge T \right)^{\frac{1}{2}}.
$$
Note that we can get a
similar control on $\int_X df \wedge d^c f \wedge \om_{\p} \wedge T$.
Thus $(\dag)_T$ follows from the following last observation
$$
\int_X df \wedge d^c f \wedge S=\frac{1}{4} \int_X (\p-\f)(\om_{\f}-\om_{\p}) \wedge S
\leq \int_X (-h) \om_h \wedge S,
$$
where $S$ is any positive closed current of bidimension $(1,1)$.

It remains to treat the case where $\p \in {\mathcal E}^1(X,\om)$
but $\f$ a priori merely belongs to ${\mathcal E}(X,\om)$.
We normalize $\f,\p$ by requiring $\sup_X \f=\sup_X \p$,
hence we need to show that $\f \equiv \p$.
Set $\f_j:=\max(\f,\p-j)$. It follows from Lemma 2.3 that
$\f_j \in {\mathcal E}^1(X,\om)$. By Proposition 3.4 below,
$(\om+dd^c \f_j)^n=(\om+dd^c \p)^n$.
Now $\sup_X \f_j=\sup_X \f=\sup_X \p$ if $j \geq 0$, hence
$\f_j \equiv \p$ by previous analysis.
Letting $j \rightarrow +\infty$, we infer $\f \equiv \p$.
\end{proof}

Regarding uniqueness in the class ${\mathcal E}(X,\om)$, we have
the following observation, of independent interest.

\begin{pro}
Assume $\f,\p \in {\mathcal E}(X,\om)$ are such that $(\om+dd^c \f)^n=(\om+dd^c \p)^n$.
Then 
$$
(\om+dd^c \max[\f,\p])^n=(\om+dd^c \f)^n=(\om+dd^c \p)^n.
$$
\end{pro}

\begin{proof}
Applying Corollary 1.10 with $\mu=\om_{\f}^n=\om_{\p}^n$ yields
$$
(\om+dd^c \max[\f,\p])^n \geq \mu=(\om+dd^c {\f})^n,
$$
whence equality, since these measures have the same mass.
\end{proof}

\subsection{Quasi-homogeneous weights}

In this section we fix a weight $\chi \in {\mathcal W}_M^+$, $M \geq 1$.
Many results in this section are inspired
by their local analogues, obtained by U.Cegrell [Ce 1] when
$\chi(t)=-(-t)^p$, $p \geq 1$.

The following result is the ${\mathcal W}_M^+$-version of Lemma 2.3.

\begin{lem}
Let $\f,\p \in PSH(X,\om) \cap L^{\infty}(X)$ 
with $\f \leq \p \leq 0$. Then
$$
0 \leq \int_X (-\chi) \circ \p (\om+dd^c \p)^n \leq 
(M+1)^n \int_X (-\chi) \circ \f (\om+dd^c \f)^n.
$$
\end{lem}

\begin{proof}
The proof follows from a repeated application of the following
inequality: if $T$ is a positive closed current
of bidimension $(1,1)$, then
\begin{equation}
0 \leq \int_X (-\chi) \circ \f \,\om_{\p} \wedge T
\leq (M+1) \int_X (-\chi) \circ \f \,   \om_{\f} \wedge T.
\end{equation}

Integrating by parts yields
$$
\int_X (-\chi) \circ \f \,\om_{\p} \wedge T
=\int_X (-\chi) \circ \f \,\om \wedge T
+\int_X (-\p) dd^c(\chi \circ \f)  \wedge T.
$$
The first integral in the RHS is bounded from above 
by $\int_X (-\chi) \circ \f \,   \om_{\f} \wedge T$ because
$\int_X \chi' \circ \f \, d \f \wedge d^c \f \wedge T \geq 0$.
Observe now that
$$
dd^c \chi \circ \f =\chi'' \circ \f d\f \wedge d^c \f+\chi' \circ \f dd^c \f
\leq \chi' \circ \f \, \om_{\f},
$$
thus
$$
\int_X (-\p) dd^c(\chi \circ \f)  \wedge T
\leq \int_X (-\p) \chi' \circ \f \, \om_{\f} \wedge T
\leq M \int_X (-\chi) \circ \f \,\om_{\f} \wedge T,
$$
since $(-\p) \chi' \circ \f \leq (-\f) \chi' \circ \f \le M (-\chi) \circ \f$.
This yields (11).
\end{proof}

We let the reader check that Corollary 2.4, Theorem 2.6 and Corollary 2.7
hold when $\chi \in {\mathcal W}_M^+$,
with exactly the same proof.
We now establish the 
high-energy version of Proposition 2.5.

\begin{pro}
Let $T$ be a positive closed current of bidimension (j,j) on $X$, $0 \leq j \leq n$,
and let $\f,\p \leq 0$ be bounded $\om$-psh functions. Then
$$
0 \leq \int_X (-\chi) \circ \f \, \om_{\p}^j \wedge T
\leq 2M  \int_X (-\chi) \circ \f \, \om_{\f}^j \wedge T+
2M \int_X (-\chi) \circ \p \, \om_{\p}^j \wedge T.
$$
\end{pro}

\begin{proof}
The proof is the same as that of Proposition 2.5 except that we have to
replace inequality $0 \leq \chi'(2t) \leq \chi'(t)$ by
$0 \leq \chi'(2t) \leq M  \chi'(t)$.
The latter follows from the concavity property of $\chi$ which 
(together with $\chi(0)=0$) yields
$|\chi(t)| \leq |t| \chi'(t)$ 
and $|\chi(2t)| \leq 2|\chi(t)|$ for all $t \in \R^-$.
Therefore
$$
0 \leq \frac{\chi'(2t)}{\chi'(t)}=
\left| \frac{2t \chi'(2t)}{\chi(2t)} \right| \cdot \left| \frac{\chi(2t)}{2 \chi(t)} \right|
\cdot \left| \frac{\chi(t)}{t \chi'(t)} \right| \leq M.
$$
\end{proof}

Our next observation is useful to establish
convexity properties of the class $\E$.

\begin{lem}
For all $0 \leq \e \leq 1$ and for all $t<-1$,
$$
0 \leq \e^M |\chi(t)| \leq |\chi(\e t)| \leq \e |\chi(t)|.
$$
\end{lem}

The proof follows easily from the concavity of $\chi$,
the normalization $\chi(0)=0$ and the definition of ${\mathcal W}_M^+$.
These inequalities can be interpreted as a 
weak-homogeneity property satisfied by the weights
$\chi \in {\mathcal W}^+_M$. This allows
us to show that $\E$ is always convex in this case.

\begin{pro}
Fix $M>0$ and $\chi \in {\mathcal W}^+_M$.
There exists $C_{\chi}>0$ such that for all 
$0 \geq  \f_0,\ldots,\f_n \in PSH(X,\om) \cap L^{\infty}(X)$,
$$
0 \leq \int_X (-\chi)\circ \f_0 \, \om_{\f_1} \wedge \cdots \wedge \om_{\f_n}
\leq C_{\chi} \max_{0 \leq j \leq n} \left[ \int_X (-\chi) \circ \f_j \,  \om_{\f_j}^n \right].
$$

In particular the class $\E$ is starshaped and convex.
\end{pro}

\begin{proof}
The proposition follows from the following inequality,
\begin{equation}
0 \leq \int_X (-\chi)\circ \f_0 \, \om_{\f_1} \wedge \cdots \wedge \om_{\f_n} \leq
\frac{4nM}{\e^n [1-2n \e M]} \max_{0 \leq j \leq n}  E_{\chi}(\f_j).
\end{equation}
The proof is identical to that of Proposition 2.10
except that Proposition 2.5 has to be replaced by Proposition 3.6,
and the subadditivity and homogeneity of
$x \in \R^+ \mapsto x^p \in \R^+$, $0 \leq p \leq 1$, has to
be replaced by the convexity property of
$x \in \R^+ \mapsto (-\chi)(-x) \in \R^+$, which yields
$$
E_{\chi}\left(\e \sum_{i=1}^n \f_i\right) \leq \e \sum_{i=1}^n \int_X (-\chi) \circ \f_i \om_u^n,
$$
for $u:=\e \sum_{i=1}^n$ and 
$0 <\e<[2nM_{\chi}]^{-1}$.
\end{proof}

We now give a high-energy version of Lemma 2.11.
Note that our proof of Lemma 2.11 uses the full homogeneity of $x \mapsto x^p$.
Let us say that a weight $\chi \in {\mathcal W}^+$
is {\bf quasi-homogeneous} if there exists $C,M \geq 1$ and $0 \leq q<1$ such that
for all $0 \leq \e \leq 1$ and for all $t \leq -1$,
$$
0 \leq C^{-1} \e^{M} |\chi(t)| \leq |\chi(\e t)| \leq C \e^{M-q} [\chi(t)|.
$$

The functions $\chi_p(t)=-(-t)^p$, $p \geq 1$, belong to ${\mathcal W}^+$
and obviously satisfy the previous condition.
Here again we shall use  the notation
$$
{\mathcal E}^p(X,\om):=\E,
\text{ when } \chi(t)=-(-t)^p.
$$
These classes have been studied in a local context by U.Cegrell [Ce 1].
They are easier to understand thanks to
the homogeneity property of the weight.

\begin{lem}
Let $\chi \in {\mathcal W}_M^+$ be a quasi-homogeneous weight and
let $\mu$ be a probability measure on $X$. 

Then $\chi \circ \E \subset L^1(\mu)$
if and only if there exists $C>0$ such that for all functions $\f \in PSH(X,\om) \cap L^{\infty}(X)$
normalized by $\sup_X \f=-1$, one has
$$
0 \leq \int_X (-\chi) \circ \f \,  d\mu \leq C 
\left( \int_X (-\chi) \circ\f \, \om_{\f}^n\right)^{\g},
$$
where $0<\g:=M/(M-q+1)<1$.
\end{lem}

The proof, in the same vein as that of Lemma 2.11, is left to the reader.

\begin{exa}
The weight $\chi(t)=-(-t)^p [\log(e-t)]^a\in {\mathcal W}^+$,$a>0$, is an example of
quasihomogeneous weight which is not homogeneous.
\end{exa}

\section{Range of the complex Monge-Amp\`ere operator}

We now turn to the central question of describing the range of the
complex Monge-Amp\`ere operator on classes $\E$, and on class
${\mathcal E}(X,\om)$.

An obvious necessary condition for solving the Monge-Amp\`ere equation
$$
(MA)_{\mu} \hskip2cm (\om+dd^c \f)^n= \mu,  
$$
is that
$\mu$ should be a positive Radon measure of total mass $\int_X \om^n$ on $X$.
For simplicity we assume throughout the rest of this section that
$\om$ is a {\it K\"ahler form}, normalized by $\int_X \om^n=1$, 
and that $\mu$ is a probability measure.

\subsection{The classes ${\mathcal E}_{\chi}(X,\om)$}

We fix here an increasing function $\chi:\R^- \rightarrow \R^-$
such that $\chi \in {\mathcal W}^- \cup  {\mathcal W}^+_M$, $M \geq 1$.
It follows from Propositions 2.5 and 3.6 that if $\mu=(\om+dd^c \f)^n$
for some function $\f \in \E$, then
$$
\! \! \!  \! \! \!  \! \! \!  \! \! \!  \! \! \! 
\! \! \!  \! \! \!  \! \! \!  \! \! \!  \! \! \! 
\! \! \!  \! \! \!  \! \! \!  \! \! \!  \! \! \! 
\! \! \!  \! \! \!  \! \! \!  \! \! \!  \! \! \! 
\! \! \!  \! \! \!  \! \! \!  \! \! \!  \! \! \! 
(I_{\chi}) \hskip3cm \chi \circ \E \subset L^1(\mu).
$$
We now show that the converse is true under a quantitative
version of $(I_{\chi})$.

\begin{thm}
Suppose there exists $F:\R^+\rightarrow \R^+$ with $\lim_{+\infty} F(t)/t=0$,
such that for all $\p \in PSH(X,\om) \cap L^{\infty}(X)$, $\sup_X \p=-1$,
$$
\! \! \!  \! \! \!  \! \! \!  \! \! \!  \! \! \! 
\! \! \!  \! \! \!  \! \! \!  \! \! \!  \! \! \! 
\! \! \!  \! \! \!  \! \! \!  \! \! \!  \! \! \! 
\! \! \!  \! \! \!  \! \! \!  \! \! \! \! 
(II_{\chi}) \hskip2cm
0 \leq \int_X (-\chi) \circ \p  d\mu \leq F\left( E_{\chi}(\p) \right).
$$

Then there exists $\f \in \E$ such that
$$
\mu=(\om+dd^c \f)^n \; \text{ and } \; \sup_X \f=0.
$$
\end{thm}

It is an interesting problem to determine whether conditions
$(I_{\chi})$ and $(II_{\chi})$ are equivalent (obviously
$(II_{\chi})$ implies $(I_{\chi})$).
This is the case  when
$\chi(t)=-(-t)^p$, $p>0$, thanks to Lemmas 2.11, 3.9,
with $F(t)=t^{p/(p+1)}$.
Thus we obtain a complete characterization of the range
of the Monge-Amp\`ere operator in this case
(Theorem C in our introduction).

\begin{thm}
Let $\mu$ be a probability measure on $X$ and $p>0$.
There exists $\f \in {\mathcal E}^p(X,\om)$ such that
$\mu=(\om+dd^c \f)^n$ 
if and only if ${\mathcal E}^p(X,\om) \subset L^p(\mu)$.
\end{thm}

The proof of Theorem 4.1 will occupy the rest of section 4.1.
The proof follows the lines of U.Cegrell's one, in the local case [Ce 1]:
\begin{itemize}

\item We approximate $\mu$ by smooth probability volume forms $\mu_j$
 using local convolutions and a partition of unity.

\item We invoke Yau's solution of the Calabi conjecture to find
 uniquely determined $\om$-psh functions $\f_j$ such that
 $\mu_j=\om_{\f_j}^n$, $\sup_X \f_j=-1$.

\item Since $\om$-psh functions $\f$ normalized by $\sup_X \f=-1$ form a compact subset
 of $L^1(X)$, we can assume that $\f_j \rightarrow \f$ in $L^1(X)$.

\item The quantitative integrability condition $(II_{\chi})$ guarantees,
at least when $\chi(t)=t$, 
$\sup_j \int (-\chi) \circ \f_j \,\om_{\f_j}^n <+\infty$,
 hence yields $\f \in {\mathcal E}_{\chi}(X,\om)$.

\item The delicate point is then to show that $\om_{\f_j}^n \rightarrow \om_{\f}^n$. This is done
 by showing that $\int |\f_j-\f| d\mu_j \rightarrow 0$ and invoking Proposition 2.8.
 We need here to assume first that $\mu$ is suitably dominated by the 
 Monge-Amp\`ere capacity (in the spirit of [K 1]).

\item We then treat the general case by  using a Radon-Nikodym decomposition of the measure $\mu$.
\end{itemize}

Here follow the technical details. 
Let $\mu$ be an arbitrary probability measure on $X$.
Let $\{U_i\}$ be a finite covering of $X$ by open sets 
$U_i$ which are biholomorphic to the unit ball of $\C^n$. In each $U_i$ we let
$\mu_{\e}^{U_i}:=\mu_{|U_i} * \r_{\e}$ denote local regularization of $\mu_{|U_i}$ by means
of convolution with radial nonnegative smooth approximations $\r_{\e}$ of the Dirac mass.
Let $\{\theta_i\}$ be a partition of unity subordinate to $\{U_i\}$ and set
$$
\mu_j:=c_j \left[ \sum_i \theta_i \mu_{\e_j}^{U_i}+\e_j \om^n \right],
$$
where $\e_j \searrow 0$ and $c_j \nearrow 1$ is chosen so that $\mu_j(X)=1$.
Thus the $\mu_j$'s are smooth probability volume forms which converge weakly towards $\mu$.
It follows from the solution of the Calabi conjecture [Y], that there exists a unique
function $\f_j \in PSH(X,\om) \cap {\mathcal C}^{\infty}(X)$ such that
$$
\mu_j=\om_{\f_j}^n \; \; \text{ and } \; \; \sup_X \f_j=-1.
$$
Recall from Proposition 1.7 in [GZ 1] that 
${\mathcal F}:=\{ \f \in PSH(X,\om) \, / \, \sup_X \f=-1\}$
is a compact subset of $L^1(X)$. Passing to a subsequence if necessary, we can therefore assume
$\f_j \rightarrow \f$ in $L^1(X)$, where $\f \in PSH(X,\om)$ with $\sup_X \f=-1$.

\begin{lem}
There exists $C>1$ such that for all $j \in \N$,
$$
\int_X (-\f_j) \, \om_{\f_j}^n \leq C \int_X (-\f_j) \, d\mu.
$$
If ${\mathcal E}^1(X,\om) \subset L^1(\mu)$, then
$\sup_j \int_X (-\f_j) \, \om_{\f_j}^n <+\infty$, hence
$\f \in {\mathcal E}^1(X,\om)$.
\end{lem}

\begin{proof}
Since $c_j \rightarrow 1$ and $\e_j \rightarrow 0$, we can write
$$
\int_X (-\f_j) \, \om_{\f_j}^n=
\sum_i \int_X \theta_i (-\f_j) \, d\mu_{\e_j}^{U_i}+o(1),
$$
where 
$$
\int_X \theta_i (-\f_j) \, d \mu_{\e_j}^{U_i} \leq 
\int_{U_i} (-\f_j*\r_{\e_j}) d\mu.
$$
Now $\f_j=u_j^i-\g_i$ in $U_i$, where $\g_i$ is a smooth local potential of $\om$ in $U_i$ and
$u_j^i$ is psh in $U_i$. Therefore
$-u_j^i*\r_{\e_j} \leq -u_j^i$, while $\g_i * \r_{\e_j}$ converges 
uniformly towards $\g_i$. We infer
$$
\int_{U_i} (-\f_j*\r_{\e_j}) d\mu \leq 
\int_{U_i} (-\f_j)  d\mu +o(1)
$$
hence 
$$
\int (-\f_j) \, d\mu_j= 
\int (-\f_j) \, \om_{\f_j}^n \leq C \int (-\f_j) \,d\mu.
$$

When ${\mathcal E}^1(X,\om) \subset L^1(\mu)$, it follows from Lemma 2.11
(case $p=1$) that 
$$
\int_X (-\f_j) \, \om_{\f_j}^n  \leq C' \left( \int_X (-\f_j) \, \om_{\f_j}^n  \right)^{1/2},
$$
hence $\sup_j \int_X (-\f_j) \, \om_{\f_j}^n \leq (C')^2<+\infty.$
Thus $\f \in {\mathcal E}^1(X,\om)$.
\end{proof}

We now would like to apply Proposition 2.8
to insure that $\mu=\om_{\f}^n$. For this we 
first assume that $\mu$ belongs to the compact convex set
${\mathcal M}_A$ of probability measures $\nu$ on $X$
which satisfy
$$
\nu(K) \leq A Cap_{\om}(K), \; \text{ for all Borel set } K \subset X.
$$
Here $A>0$ is a fixed constant.

\begin{lem}
Assume $\mu \in {\mathcal M}_A$, then
$$
\int_X \f_j d\mu \rightarrow \int_X \f d\mu \; \; 
\text{ and } \; \; 
\int_X |\f_j-\f|d\mu_j \rightarrow 0.
$$
\end{lem}

\begin{proof}
When the $\f_j$'s are uniformly bounded, the first convergence follows 
from standard arguments (see [Ce 1]). Set
$$
\f_j^{(k)}:=\max(\f_j,-k) \; \; \text{ and } \; \; \f^{(k)}:=\max(\f,-k).
$$
We will be done with the first convergence if we can show that 
$\int |\f_j^{(k)}-\f_j|d\mu \rightarrow 0$
uniformly in $j$ as $k \rightarrow +\infty$. This is where we use 
our assumption on $\mu$, since 
$$
\int_X  |\f_j^{(k)}-\f_j|d\mu \leq 
2 \int_{(\f_j<-k)} (-\f_j)d\mu \leq \frac{2}{\sqrt{k}} 
\int_X (-\f_j)^{3/2} \, d\mu \leq \frac{C}{\sqrt{k}}.
$$
Indeed Proposition 5.3 shows that ${\mathcal E}^1(X,\om) \subset L^1(\mu)$,
hence $\sup_j \int_X (-\f_j) \, \om_{\f_j}^n <+\infty$ by previous Lemma.
Moreover
$$
\int_X (-\f_j)^{3/2} d\mu=1+\frac{3}{2} \int_1^{+\infty} \sqrt{t} \mu(\f_j<-t) dt
\leq 1+\frac{3A}{2} \int_1^{+\infty} \sqrt{t}Cap_{\om}(\f_j<-t) dt.
$$
It follows now from Lemma 5.1 that
$$
Cap_{\om}(\f_j<-t) \leq C' t^{-2} \int_X (-\f_j) \om_{\f_j}^n \leq C'' t^{-2},
$$
which proves that $\int_X (-\f_j)^{3/2} d\mu$ is uniformly bounded from above.

It remains to prove a similar convergence when $\mu $ is replaced by $\mu_j$.
It actually suffices to consider the case of measures
$\mu_j^U:=\mu_{|U}*\r_{\e_j}$. Now
$$
\int_U |\f_j-\f|d\mu_j^U=\int_U \left(\int_U |u_j(\zeta)-u(\zeta)| 
\r_{\e_j}(z-\zeta) d\lambda(\zeta) \right) d\mu(z),
$$
where as above, $u_j$, $u$ are psh functions in $U$ such that $\f_j=u_j-\g$ and  $\f=u-\g$ in $U$,
$\g$ is a local potential of $\om$ in $U$ and $d\lambda$ denotes the Lebesgue measure in $U$.
The lemma will be proved if we can show that $\int w_j d\mu \rightarrow 0$, where
$$
w_j(z):=\int_U |u_j(\zeta)-u(\zeta)| \r_{\e_j}(z-\zeta) d\lambda(\zeta).
$$
Define $\tilde{u}_j:=(\sup_{k \geq j} u_k )^*$. This is a sequence of psh functions in $U$ which decrease
towards $u$. Observe that $\tilde{u}_j \geq \max(u,u_j)$ so that
$$
w_j \leq 2 \tilde{u}_j * \r_{\e_j} -u * \r_{\e_j}-u_j* \r_{\e_j}
\leq 2(\tilde{u}_j*\r_{\e_j}-u)+(\f-\f_j).
$$
It follows from the monotone convergence theorem that
$\int (\tilde{u}_j*\r_{\e_j}-u) d\mu \rightarrow 0$, while 
$\int (\f_j-\f)d\mu \rightarrow 0$ by the first part of lemma.
Therefore $\int w_j d\mu \rightarrow 0$ and we are done.
\end{proof}

What we have shown so far is that for any measure $\nu \in {\mathcal M}_A$,
$A>0$, there exists a unique $\f \in {\mathcal E}^1(X,\om)$ such that
$\nu=(\om+dd^c \f)^n$ and $\sup_X \f=-1$.

We now proceed with the proof of Theorem 4.1. For this we need the following
observation:

\begin{lem}
Let $\mu$ be a probability measure which does not charge pluripolar sets.
Then there exists $u \in PSH(X,\om) \cap L^{\infty}(X)$
and $0 \leq f \in L^1(\om_u^n)$ such that
$\mu=f (\om+dd^c u)^n$.
\end{lem}

\begin{proof}
Recall that ${\mathcal M}_1$
is a compact convex subset of the set of all
probability measures on $X$.
Let $\mu$ be a probability measure which does not charge pluripolar sets.
It follows from a generalization of Radon-Nikodym theorem [R] that
$$
\mu=f_1 \nu+\sigma,
\text{ where } \nu \in {\mathcal M}_1, \;
0 \leq f_1 \in L^1(\nu), \, \text{ and }
\sigma \perp {\mathcal M}_1.
$$
Now $\sigma$ is carried by a pluripolar set since 
${\mathcal M}_1$ contains all measures
$\om_u^n$, $0 \leq u \leq 1$, hence $\sigma=0$.
By our previous analysis, $\nu=(\om+dd^c v)^n$,
where $v \in {\mathcal E}^1(X,\om)$, $\sup_X v=0$.
Set $u:=\exp v$. This is again a $\om$-psh function since
$$
\om_u=e^v \om_v +[1-e^v] \om+e^v dv \wedge d^c v \geq e^v \om_v \geq 0.
$$
Observe that $0 \leq u \leq 1$ and $\om_u^n \geq e^{nv} \om_v^n$, thus
$\mu$ is absolutely continuous with respect to $\om_u^n$.
\end{proof}

We now go on with  the proof of Theorem 4.1. 
Let $\mu$ be a probability measure which satisfies $(II)_{\chi}$.
Then $\mu$ does not charge pluripolar sets, hence it
writes $\mu=f \om_u^n$, where 
$u \in PSH(X,\om) \cap L^{\infty}(X)$ is so that $0\leq u \leq 1$.
Consider
$$
\mu_j:= \d_j \min(f,j) (\om+dd^c u)^n,
$$
where $\d_j \searrow 1$ so that $\mu_j(X)=1$. 
Note that
$$
\mu_j \leq j \d_j Cap_{\om}, 
\text{ since } 0 \leq u \leq 1,
$$
thus $\mu_j \in {\mathcal M}_{j \d_j}$. 
It follows therefore from previous analysis that
there exists a unique $\f_j \in {\mathcal E}^1(X,\om)$ with $\sup_X \f_j =-1$ 
and $\mu_j=\om_{\f_j}^n$.
We can assume $\f_j \rightarrow \f$ in $L^1(X)$
and $\d_j \leq 2$.

Suppose  first that $\chi \in {\mathcal W}^-$, so that
$\f_j \in {\mathcal E}_{\chi}(X,\om) \supset {\mathcal E}^1(X,\om)$.
It follows from Theorem 2.6 that $(II_{\chi})$ can be applied to
the functions $\f_j$, hence
$$
\int (-\chi) \circ \f_j \,  \om_{\f_j}^n \leq \d_j \int (-\chi \circ \f_j) d\mu 
\leq 2 C F \left(\int_X  (-\chi) \circ \f_j \,  \om_{\f_j}^n \right).
$$
Thus $\sup_j E_{\chi}(\f_j)<+\infty$, hence $\f \in {\mathcal E}_{\chi}(X,\om)$
by Corollary 2.7.
We set
$$
\Phi_j:= (\sup_{k \geq j} \f_k)^* \; \; \text{ and } F_j:=\inf_{k \geq j} \d_k \min(f,k).
$$
Clearly $\Phi_j \in {\mathcal E}_{\chi}(X,\om)$ with 
$\Phi_j \searrow \f$ and $F_j \nearrow f$. It follows 
therefore from 
Corollary 1.10 that
$$
\om_{\Phi_j}^n \geq F_j \nu.
$$
We infer $\om_{\f}^n \geq \mu$, whence equality since these are both probability measures.
Thus the proof of Theorem 4.1 is complete when $\chi \in {\mathcal W}^-$.

Assume now $\chi \in {\mathcal W}^+_M$. We can apply previous reasoning
if we can show that $\f_j \in {\mathcal E}_{\chi}(X,\om)$ for all $j$.
This is the case, as we claim that the $\f_j$'s actually belong
to ${\mathcal E}^p(X,\om)$ for all $p \geq 1$. Indeed recall that
$$
\mu_j=(\om+dd^c \f_j)^n \leq j \d_j(\om+dd^c u)^n \leq A_j Cap_{\om},
$$
hence $\mu_j \in {\mathcal M}_{A_j}$. Observe that for all $p \geq 1$,
$$
(*) \hskip2cm
{\mathcal E}^p(X,\om) \subset L^{p+1/2}(\nu),
\text{ for all  } \nu \in {\mathcal M}_A.
$$
This has been established when $p=1$ in the proof of Lemma 4.4.
Since $\f_j \in {\mathcal E}^1(X,\om)$
and $\mu_j=(\om+dd^c \f_j)^n \in {\mathcal M}_{A_j}$, we infer
$\f_j \in L^{3/2}([\om+dd^c \f_j]^n)$, i.e. 
$\f_j \in {\mathcal E}^{3/2}(X,\om)$.
We can thus use $(*)$ inductively to obtain $\f_j \in \cap_p {\mathcal E}^p(X,\om)$.
The proof of Theorems 4.1 and 4.2 is now complete.

\subsection{Non pluripolar measures}

We now describe the range of the complex Monge-Amp\`ere
operator on the class ${\mathcal E}(X,\om)$ (Theorem A).

\begin{thm}
There exists
$\f \in {\mathcal E}(X,\om)$ such that $\mu=(\om+dd^c \f)^n$ if and only if
$\mu$ does not charge pluripolar sets.
\end{thm}

\begin{proof}
One implication is obvious. Namely if $\mu=(\om+dd^c \f)^n$ for some function
$\f \in {\mathcal E}(X,\om)$, then  $\mu$ does not charge pluripolar sets,
as follows from Theorem 1.3.
\vskip.1cm

Assume now that $\mu$ does not charge pluripolar sets.
It follows from Lemma 4.5 that we can find 
$u \in PSH(X,\om) \cap L^{\infty}(X)$
and $0 \leq f \in L^1(\om_u^n)$ such that
$\mu=f \om_u^n$. Set
$$
\mu_j:=c_j \min(f,j) \om_u^n,
$$
where $c_j \searrow 1$ is such that $\mu_j(X)=\mu(X)=1$. 
We can assume without loss of generality that $1 \leq c_j \leq 2$.

It follows from Chern-Levine-Nirenberg inequalities (see Proposition 3.1 in [GZ 1])
that ${\mathcal E}^1(X,\om) \subset PSH(X,\om) \subset L^1(\mu_j)$,
because $\mu_j\leq 2j \om_u^n$ and $u$ is bounded.
Thus by Theorem 4.2, case $p=1$, there exists a unique function $\f_j \in {\mathcal E}^1(X,\om)$
such that
$$
\mu_j=(\om+dd^c \f_j)^n \; \text{ and } \; \sup_X \f_j=0.
$$
Passing to a subsequence if necessary, we can assume
$\f_j \rightarrow \f$ in $L^1(X)$, for some function
$\f \in PSH(X,\om)$ such that $\sup_X \f=0$.
\vskip.1cm

We claim that $\f \in {\mathcal E}_{\chi}(X,\om)$ for some
$\chi \in {\mathcal W}^-$,
and that $\mu=(\om+dd^c {\f})^n$.
The function $\chi$ is defined as follows.
Let $\g:\R^+ \rightarrow \R^+$ be a convex increasing function with
$\lim_{z \rightarrow +\infty} \g(z)/z=+\infty$ such that
$\g \circ f$ still belongs to $L^1(\om_u^n)$
(see [RR] for the construction of $\g$).
Let $\g^*$ be the Young-conjugate function of $\g$,
$$
\g^*:z \in \R^+ \mapsto \sup \{ zy-\g(y) \, / \, y \in \R^+\} \in \R^+,
$$
and set $\chi(t):=-(\g^*)^{-1}(-t)$: this is a convex increasing function such that
$\chi(-\infty)=-\infty$, which satisfies
\begin{equation}
(-\chi)(-t) \cdot  f(x) \leq -t+\g \circ f(x),
\;
\text{ for all } (t,x) \in \R^- \times X,
\end{equation}
since $z y \leq \g^*(z)+\g(y)$ for all $z,y \in \R^+$.

Note that we can assume $\chi(0)=0$ without loss of generality,
by imposing $\g(0)=0$.
Applying (13) with $t=\f_j(x)$ and averaging against $\om_u^n$,
we infer
$$
0 \leq \int_X (-\chi)\circ \f_j \, \om_{\f_j}^n \leq
2 \int_X (-\chi) \circ \f_j \, f \om_u^n
\leq 2 \int_X (-\f_j) \om_u^n+2\int_X \g \circ f \, \om_u^n.
$$

Again it follows from Chern-Levine-Nirenberg inequalities 
that $\int_X (-\f_j) \om_u^n$ is uniformly bounded from above,
because $u$ and $\int_X (-\f_j) \om^n$ are bounded.
Thus the sequence of integrals
$(\int_X (-\chi) \circ \f_j \, \om_{\f_j}^n )$ is bounded.

Set $\phi_j:=(\sup_{l \geq j} \f_l)^* \in {\mathcal E}^1(X,\om)$.
It follows from Corollary 1.10 that
$$
\om_{\phi_j}^n \geq \inf_{l \geq j} \om_{\f_l}^n \geq \min(f,j) \om_u^n.
$$
The sequence $(\phi_j)$ decreases towards $\f$ and satisfies
$\f_j \leq \phi_j \leq 0$, hence by the fundamental inequality,
$$
0 \leq \int_X (-\chi) \circ \phi_j \, \om_{\phi_j}^n
\leq 2^n \int_X (-\chi) \circ \f_j \, \om_{\f_j}^n
\leq M_0<+\infty.
$$
Therefore $\f\in \E$ and
$$
\om_{\f}^n=\lim_{j \rightarrow +\infty} \om_{\phi_j}^n
\geq \lim_{j \rightarrow +\infty} \min(f,j) \om_u^n=\mu.
$$
Since $\mu$ and $\om_{\f}^n$ are both probability measures, 
this actually yields equality,
hence the proof is complete.
\end{proof}

\begin{rqe}
Note that when $\mu=f \om^n$ is a measure with density $0 \leq f$
such that $f \log f \in L^1(\om^n)$,
our proof shows that $\mu=(\om+dd^c \f)^n$ with  $\f \in {\mathcal E}_{\chi}(X,\om)$,
for $\chi(t)=-\log(1-t)$
(this is a critical case in the Orlicz classes considered by S.Kolodziej [K 1]).
\end{rqe}

\section{Examples and Applications}

\subsection{Capacity of sublevel sets}

Recall that the Monge-Amp\`ere capacity associated to $\om$ is
defined by
$$
Cap_{\om}(K):=\sup \left\{ \int_K \om_{u}^n \; / \, u \in PSH(X,\om),
\, 0 \leq u \leq 1 \right\},
$$
where $K$ is any Borel subset of $X$. This capacity vanishes on pluripolar
sets, more precisely $Cap_{\om}(\f <-t) \leq C_{\f}/t$ for every fixed $\om$-psh
function $\f$. This estimate is sharp in the sense that
$Cap_{\om}(\f <-t) \geq C'_{\f}/t$ when $\om_{\f}$ is the current of integration
along a hypersurface.
However when $\f$ belongs to ${\mathcal E}_{\chi}(X,\om)$, one can establish finer
estimates as our next result shows.

\begin{lem}
Fix $\chi \in {\mathcal W}^- \cup {\mathcal W}_M^+$, $M \geq 1$.
If $\f \in {\mathcal E}_{\chi}(X,\om)$,
then 
$$
\exists C_{\f}>0, \forall t>1, \; 
Cap_{\om}(\f <-t) \leq C_{\f} |t \, \chi(-t)|^{-1}.
$$

Conversely if there exists $C_{\f}, \e>0$ such that 
for all $t>1$,
$$
Cap_{\om}(\f <-t) \leq C_{\f} | t^{n+\e} \, \chi(-t)|^{-1},
$$ 
then $\f \in {\mathcal E}_{\chi}(X,\om)$.
\end{lem}

The proof relies on the comparison principle, as it was used in different
contexts by U.Cegrell and S.Kolodziej (see [Ce 1], [K 1], [CKZ]).

\begin{proof}
Fix $\f \in \E$, $\f \leq -1$ and $u \in PSH(X,\om)$
with $-1 \leq u \leq 0$.
For $t \geq 1$, observe that $\f/t \in {\mathcal E}(X,\om)$ and
$
(\f <-2t) \subset (\f/t<u-1) \subset (\f<-t).
$
It therefore follows  from the comparison principle that 
\begin{eqnarray*}
\lefteqn{\! \! \! \! \! \! \! \! \! \! 
\! \! \! \! \!  \! \! \! \! \! 
\! \! \! \! \!  \! \! \! \! \! 
\int_{(\f<-2t)} (\om+dd^c u)^n \leq \int_{(\f<-t)} (\om+t^{-1}dd^c \f)^n } \\
&\leq& \int_{(\f<-t)} \om^n +t^{-1} \sum_{j=1}^n \left( \begin{array}{cc} n \\ j \end{array} \right)
\int_{(\f<-t)} \om_{\f}^j \wedge \om^{n-j}.
\end{eqnarray*}
Recall  now that $Vol_{\om}(\f<-t)$ decreases exponentially
fast [GZ 1], and observe that for all $1 \leq j \leq n$,
$$
\int_{(\f<-t)} \om_{\f}^j \wedge \om^{n-j}
\leq \frac{1}{|\chi(-t)|} \int_X (-\chi) \circ \f \, \om_{\f}^j \wedge \om^{n-j}
\leq \frac{1}{|\chi(-t)|} E_{\chi}(\f).
$$
This yields our first assertion.

The second assertion follows from similar considerations. Namely for
$\f \in PSH(X,\om)$, $\f \leq -1$, set $\f_t:=\max(\f,-t)$.
Then $u=\f_t/t \in PSH(X,\om)$ with $-1 \leq u \leq 0$ for all $t \geq 1$.
Since $\om_u^n \geq t^{-n} \om_{\f_t}^n$, we infer
$$
t^{-n} \om_{\f_t}^n(\f<-t) \leq \om_u^n(\f<-t) \leq Cap_{\om}(\f<-t).
$$
If $Cap_{\om}(\f<-t) \leq C t^{-n-\e} |\chi(-t)|^{-1}$, this
shows that $\om_{\f_t}(\f<-t) \rightarrow 0$, hence $\f \in {\mathcal E}(X,\om)$.
Moreover
$$
\om_{\f}^n(\f \leq -t)=\int_X \om^n -\om_{\f}^n(\f>-t)
=\int_X \om^n -\om_{\f_t}^n(\f>-t)=\om_{\f_t}^n(\f \leq -t)
$$
thus $\om_{\f}^n(\f \leq -t) \leq C t^{-\e} |\chi(-t)|^{-1}$.
This yields
$$
\int_X (-\chi) \circ \f \, \om_{\f}^n
=\int_1^{+\infty} \chi'(-t) \om_{\f}^n(\f<-t) dt
\leq C_{\f} \int_1^{+\infty} \frac{t\chi'(-t)}{|\chi(-t)|} \frac{1}{t^{1+\e}} dt
<+\infty,
$$
so that $\f \in \E$.
\end{proof}

These estimates allow us to give several examples of functions which belong
to the classes ${\mathcal E}_{\chi}(X,\om)$.
We simply mention the following ones:

\begin{exa}
Let $\f$ be any $\om$-plurisubharmonic function such that $\f \leq 0$. Then 
$\p:=-\log(1-\f)=\chi \circ \f \in PSH(X,\om)$ since 
$$
dd^c \chi \circ \f=\chi" \circ \f \, d\f \wedge d^c \f+\chi' \circ \f dd^c \f
$$
and $\chi \in {\mathcal W}^-$
with $\chi'(t) \leq 1$ when $t \leq 0$.

Observe that the capacity of the sets $(\p<-t)$ decrease
exponentially fast, thus 
$\p \in \cap_{p \geq 1} {\mathcal E}^p(X,\om)$ (Lemma 5.1).
Together with Theorem 7.2 in [GZ 1], this shows that
${\mathcal E}^p(X,\om)$ characterizes pluripolar sets (for any $p \geq 1$).
\end{exa}

We now want to give some example of probability measures
that satisfy the assumptions of Theorems 4.1, 4.2, 4.6.

When $\mu=f\om^n$ has density $f \in L^r(X)$, $r>1$, S.Kolodziej has proved [K 1]
that $\mu=\om_{\p}^n$ for some {\it bounded} $\om$-psh function $\p$.
This is because $\mu$ is strongly dominated by $Cap_{\om}$ 
in this case (see Proposition 5.3 below).
When the density is only in $L^1$, this does not hold. Consider for instance
$\mu=f\om^n$, where $f \in {\mathcal C}^{\infty}(X \setminus \{a\})$ is
such that 
$$
f(z) \simeq \frac{1}{||z||^4 (-\log ||z||)^2}-1
$$
near the point $a=0$, in a local chart. Observe that
$$
\f(z):=\e \chi(z) \log ||z|| \in PSH(X,\om)
$$
if $\chi$ is a cut-off function so that $\chi \equiv 1$ near $a=0$, and $\e>0$ is small 
enough. Now $\f \notin L^1(\mu)$ but still
${\mathcal E}^1(X,\om) \subset L^1(\mu)$,
as follows from Proposition 5.3 below,
thus there exists $\p \in {\mathcal E}^1(X,\om)$ such that $\mu=\om_{\p}^n$ 
(Theorem 4.2). 

Observe also that there are measures $\mu=f \om^n$ with $L^1$-density such that
${\mathcal E}^1(X,\om) \not\subset L^1(\mu)$: one can consider for instance
$f_{\e}$ that looks locally near $a=0$ like
$[ ||z||^4 (-\log ||z||)^{1+\e}]^{-1}$, for $\e>0$ small enough.
However this measure $\mu$ does not charge pluripolar sets, hence
is a Monge-Amp\`ere measure of some function 
$\f \in \E$, $\chi \in {\mathcal W}^-$
with slower growth.

\begin{pro}
Let $\mu$ be a probability measure on $X$.

Assume there exists $\a>p/(p+1)$ and $A>0$ such that
\begin{equation}
\mu(E) \leq A Cap_{\om}(E)^{\a},
\end{equation}
for all Borel sets $E \subset X$. Then ${\mathcal E}^p(X,\om) \subset L^p(\mu)$.

Conversely assume ${\mathcal E}^p(X,\om) \subset L^p(\mu)$, $p>1$. Then there exists
$0<\a<1$ and $A>0$ such that $(14)$ is satisfied.
\end{pro}

\begin{proof}
We can assume w.l.o.g. that $\a<1$. Let $\f \in {\mathcal E}^p(X,\om)$
with $\sup_X \f=-1$. It follows from H\"older inequality that
\begin{eqnarray*}
\lefteqn{ \! \! \! \! \! \! \! \! \! 
0  \leq  \int_X (-\f)^p d\mu =1+ p \int_1^{+\infty} t^{p-1} \mu(\f<-t) dt} \\
&\leq& 1+pA \int_1^{+\infty}t^{p-1} \left[ Cap_{\om}(\f<-t) \right]^{\a} dt \\
&\leq& 1+pA \left[ \int_1^{+\infty} t^{\frac{p-\a(p+1)}{1-\a}-1} dt \right]^{1-\a}
\cdot \left[ \int_1^{+\infty} t^p Cap_{\om}(\f<-t) dt \right]^{\a}.
\end{eqnarray*}
The first integral in the last line converges since $p-\a(p+1)<0$ and 
it follows from the proof of Lemma 5.1 that the last one  is bounded from above
by $C_{p,\a} \left( \int (-\f)^p \om_{\f}^n \right)^{\a}$.
Therefore ${\mathcal E}^p(X,\om) \subset L^p(\mu)$.

Assume conversely that ${\mathcal E}^p(X,\om) \subset L^p(\mu)$.
It follows from Theorem 4.2 that $\mu=\om_{\p}^n$ for some
function $\p \in {\mathcal E}^p(X,\om)$ such that $\sup_X \p=-1$.
We claim then that there exists $\g_p \in ]0,1[$ and $A>0$ such that
for all functions $\f \in PSH(X,\om)$ with $-1 \leq \f \leq 0$, one has
\begin{equation}
0 \leq \int_X (-\f)^p \om_{\p}^n \leq A \left[ \int_X (-\f)^p \om_{\f}^n \right]^{\g_p}.
\end{equation}
We leave the proof of this claim to the reader 
(the exponent $\g_p=(1-1/p)^n$ would do).
We apply now (15) to the extremal function
$\f=h_{E,\om}^*$ introduced in [GZ 1]. It follows from
Theorem 3.2 in [GZ 1] that
$$
0 \leq \mu(E) \leq \int_X (-h_{E,\om}^*)^p d\mu \leq A  Cap_{\om}(E)^{\g_p}.
$$
\end{proof}

This proposition allows to produce several examples of measures satisfying
${\mathcal E}^p(X,\om) \subset L^p(\mu)$ as in the local theory
(see [K 1], [Z]).
It can also be used, together with Theorem 4.2, to prove that functions from the 
local classes of Cegrell ${\mathcal E}^r(\Omega)$, $\Omega$ a bounded hyperconvex
domain of $\C^n$, can be sub-extended as global functions
$\f \in {\mathcal E}^{p}(\P^n,A\om_{FS})$, for all $p<r/n$ and some $A >0$
(see [CKZ] for similar results).

\subsection{Complex dynamics}

\subsubsection{Large topological degree}
Let $f:X \rightarrow X$ be a meromorphic endomorphism whose topological degree
$d_t(f)$ is large in the sense that
$d_t(f)>r_{k-1}(f)$, where $r_{k-1}(f)$ denotes the spectral radius
of the linear action induced by $f^*$ on $H^{k-1,k-1}(X,\R)$.
In this case there exists a unique invariant measure $\mu_f$
of maximal entropy which can be decomposed as
$$
\mu_f:=\Theta+dd^c ({\mathcal T}),
$$
as was proved by the first author in [G]. Here $\Theta$
is a smooth probability measure and ${\mathcal T} \geq 0$ is
a positive current of bidegree $(k-1,k-1)$. In particular
$$
{\mathcal E}^1(X,\om) \subset PSH(X,\om) \subset L^1(\mu) ,
$$
as follows from Stokes theorem (see Theorem 2.1 in [G]
and Example 2.8 in [GZ 1]).
Here $\om$ denotes any fixed K\"ahler form on $X$, normalized 
by $\int_X \om^n=\mu(X)=1$.
It follows therefore from Theorem 4.2 and Theorem 3.3 that there
exists a unique function $g_f \in {\mathcal E}^1(X,\om)$ such that
$\sup_X g_f=-1$ and 
$$
\mu_f=(\om+dd^c g_f)^k.
$$
It is an interesting problem to establish further regularity properties
of $g_f$ in order for example to estimate the pointwise dimension of the measure $\mu_f$.
When $f:\P^n \rightarrow \P^n$ is an {\it holomorphic}
endomorphism of the complex projective space $X=\P^n$, it is known [S] that
$g_f$ is an H\"older-continuous function.

\subsubsection{Small topological degree}

Let $T$ be a positive closed current of bidegree $(p,p)$ on $X$.
Let $\om$ denote a positive closed $(1,1)$ current with bounded
potentials, and such that $\int_X T \wedge \om^{n-p}>0$.
Note that the wedge product $T \wedge \om^{n-p}$ is a well defined
positive Radon measure because $\om$ has bounded potentials [BT 3].
Fix $\chi \in {\mathcal W}^-$.

\begin{defi}
We let ${\mathcal E}_{\chi}(T,\om)$ denote the class of $\om$-psh functions
with finite $\chi$-energy with respect to $T$. This is the set of function 
$\f \in PSH(X,\om)$ for which there exists a sequence
$\f_j \in PSH(X,\om) \cap L^{\infty}(X)$  such that
$$
\f_j \searrow \f \hskip.2cm
\text{ and } \hskip.2cm
\sup_j \int_X |\chi \circ (\f_j -\sup_X \f_j)| T \wedge \left(\om+dd^c \f_j\right)^{n-p} <+\infty .
$$
\end{defi}

When $T=[X]$ is the current of integration on $X$ ($p=0$), this is the same 
notion as that of Definition 1.1.
We let the reader check that the following properties
hold, with the same proof as in the case $T=[X]$:
\begin{enumerate}
\item $\forall \f,\p \in PSH(X,\om) \cap L^{\infty}(X)$, $\f \leq \p \leq 0$ implies
$$
0 \leq \int_X (-\chi) \circ \p \,  \om_{\p}^{n-p} \wedge T
\leq 2^{n-p} \int_X (-\chi) \circ \f \,  \om_{\f}^{n-p} \wedge T;
$$
\item the operator $\om_{\f}^{n-p} \wedge T$ is well 
defined on ${\mathcal E}_{\chi}(T,\om)$
and continuous under decreasing sequences;
\item the function $\chi \circ \f$ is integrable with respect to the positive Radon measure 
$\om_{\f}^{n-p} \wedge T$.
\end{enumerate}
\vskip.3cm

This is of great interest in the dynamical study of meromorphic endomorphisms.
Assume for instance that $f:X \rightarrow X$ is a meromorphic endomorphism
of some compact k\"ahler {\it surface} $X$ ($n=2$), whose topological degree
is smaller than the first dynamical degree 
$\l=\l_1(f):=\lim_{j \rightarrow +\infty} [r_1(f^j)]^{1/j}$,
where $r_1(f)$ denotes the spectral radius of the linear action induced by $f^*$ on
$H^{1,1}(X,\R)$ (see [DF]). When $f$ is $1$-stable (i.e. when the iterates
$f^j$ of $f$ do not contract any curve into the indeterminacy locus $I_f$ of $f$),
it is known that 
$$
\frac{1}{\l^j} (f^j)^* \om^{+} \longrightarrow T^+=\om^++dd^c g^+,
$$
where $\om^+$ is a positive closed $(1,1)$-current with bounded potentials,
whose cohomology class is $f^*$-invariant, and 
$T^+$ is a canonical $f^*$-invariant positive
current, $f^* T^+=\l T^+$ (see [DDG]).
Similarly
$$
\frac{1}{\l^j} (f^j)_* \om^{-} \longrightarrow T^-=\om^-+dd^c g^-,
$$
where $\om^-$ is a positive closed $(1,1)$-current with bounded potentials,
whose cohomology class is $f_*$-invariant, and 
$T^-$ is the canonical $f_*$-invariant positive current, $f_* T^-=\l T^-$ (see [DDG]).

It is a difficult and important question to define the dynamical 
wedge product $T^+ \wedge T^-$: this product is 
expected to yield the unique measure of maximal entropy.
When $g^+ \in {\mathcal E}_{\chi}(T^-,\om^+)$, one can use our ideas above
and show that not only is the measure $\mu_f=T^+ \wedge T^-$ well defined,
but also $\chi \circ g^+ \in L^1(\mu_f)$.
In particular the measure
$\mu_f$ does not charge the set $(g^+=-\infty)$.
This has several interesting dynamical consequences, as shown
in the paper [DDG] to which we refer the reader. 
Note that when $\chi(t)=t$, this condition was introduced and
extensively studied by E.Bedford and J.Diller in [BD]: in this
case the functions $g^+,g^-$ have necessarily gradients in $L^2(X)$.
Our weaker condition allows us to handle functions whose gradient
does not necessarily belong to $L^2(X)$, as already noted above
(Example 2.14).

\vskip.2cm
To illustrate these ideas
we now consider the special case
where $X=\P^2$ is the complex projective plane 
equipped with the Fubini-Study K\"ahler form $\om$, 
and $f:\P^2 \rightarrow \P^2$ is a
$1$-stable {\it birational} map ($d_t(f)=1$), 
with $\l_1(f)>1$.  
The invariant currents write  $T^{\pm}=\om+dd^c g^{\pm}$,
where $g^{\pm} \in PSH(X,\om)$.
The functions
$g^+,g^-$ do not belong to the class ${\mathcal E}(\P^2,\om)$ because
they have positive Lelong numbers at points of indeterminacy of the mappings 
$f^n$, $n \in \Z$, however we have the following result.

\begin{pro}
If $g^+ \in {\mathcal E}^1(T^-,\om)$ then
$g:=\max (g^+,g^-) \in {\mathcal E}^1(\P^2,\om)$.
\end{pro}

\begin{proof}
We can assume without loss of generality that $g^+,g^- \leq 0$, hence $g \leq 0$.
The Bedford-Diller condition 
$g^+ \in {\mathcal E}^1(T^-,\om)$ implies that
$\nabla g^+,\nabla g^- \in L^2(\P^2)$, hence
$\nabla g \in L^2(\P^2)$ and the complex Monge-Amp\`ere
measure $(\om+dd^c g)^2$ is well defined (see [BT 2] and [BD]).
Now
$$
0 \leq \int_X (-g) \om_g^2 \leq \int_X (-g^+) \om_g^2=\int (-g^+) \om\wedge \om_g
+\int_X (-g^+) dd^c g \wedge \om_g.
$$
It follows from Cauchy-Schwarz inequality that
the first integral in the RHS is finite.
We can bound the last one from above by using Stokes theorem,
$$
\int_X (-g^+) dd^c g \wedge \om_g=\int_X (-g) dd^c g^+ \wedge \om_g
\leq \int_X (-g) \om_{g^+} \wedge \om_g.
$$
Now
\begin{eqnarray*}
\int_X (-g) \om_{g^+} \wedge \om_g 
&\leq & \int_X (-g^-) \om_{g^+} \wedge \om_g \\
&=&
\int_X (-g^-) \om_{g^+} \wedge \om+\int_X (-g^-) \om_{g^+} \wedge dd^c g \\
&\leq & O(1)+\int_X (-g) \om_{g^+} \wedge \om_{g^-} \\
&\leq & O(1)+\int_X (-g^+) \om_{g^+} \wedge T^-<+\infty,
\end{eqnarray*}
the last integral being finite because $g^+ \in {\mathcal E}^1(T^-,\om)$.
\end{proof}

It is an interesting problem to determine whether 
$\mu_f=(\om+dd^cg)^2$. This is the case when e.g. $f$ is a complex H\'enon mapping,
and it would imply -- by Theorem 4.2 -- that
${\mathcal E}^1(\P^2,\om) \subset L^1(\mu_f)$, hence 
in particular $\mu_f$ does not charge
pluripolar sets.

\subsection{Singular K\"ahler-Einstein metrics}

It is well-known that solving Monge-Amp\`ere equations
$$
\text{ }
\!\!\!\!\!\!  \!\!\! \!\!\! \!\!\! \!\!\! \!\!\! 
[MA](X,\om,\mu) \hskip2cm
(\om+dd^c \f)^n=\mu    
$$
is a way to produce K\"ahler-Einstein metrics (see [Ca], [A], [Y], [T], [K 1]).
In the classical case, the measure $\mu=f \om^n$ admits a smooth density
$f>0$. When the ambient manifold has some singularities (which is often the case
in dimension $\geq 3$), one has to allow the equation
$[MA](X,\om,\mu)$ to degenerate in two different ways:
resolving the singularities $\pi:\tilde{X} \rightarrow X$ of $X$ yields
a new equation $[MA](\tilde{X},\tilde{\om},\tilde{\mu})$, where
\begin{enumerate}
\item $\{\tilde{\om}\}=\{ \pi^* \om\}$ is a semi-positive and big class (one looses
strict positivity along the exceptional locus ${\mathcal E}_{\pi})$;
\item $\tilde{\mu}=\tilde{f} \tilde{\om}^n$ is a measure with density 
$0 \leq \tilde{f} \in L^p$, $p>1$, which may have zeroes and poles along some 
components of ${\mathcal E}_{\pi}$.
\end{enumerate}

In this paper we have focused on the second type of degeneracy. We would like to 
mention that our techniques are supple enough so that we can produce solutions
$\f \in{\mathcal E}_{\chi}(\tilde{X},\tilde{\om})$ to the Monge-Amp\`ere equations
$[MA](\tilde{X},\tilde{\om},\tilde{\mu})$, even when $\{\tilde{\om}\}$ is
merely big and semi-positive rather than K\"ahler, as was assumed for
simplicity throughout section 4. 
This is done by perturbing the form $\tilde{\om}$, adding 
small mutliples of a K\"ahler form $\Omega$,
$\om_j:=\tilde{\om}+\e_j \Omega$, and trying to understand
what happens at the limit.

An elementary but crucial observation here is that
classes $\E$ have good stability properties in the following sense.

\begin{pro}
Fix $\chi \in {\mathcal W}:={\mathcal W}^- \cup {\mathcal W}_M^+$,
$M \geq 1$.
Let $\om_j$ be a sequence of positive closed 
$(1,1)$-currents with bounded potentials, that decrease towards $\om$,
another current with bounded potentials.
Assume $\f_j \in {\mathcal E}_{\chi}(X,\om_j)$ decreases towards $\f$
pointwise and 
$\sup_j \int_X |\chi(|\f_j|)| (\om_j+dd^c \f_j)^n <+\infty$.

Then $\f \in \E$ and $(\om_j+dd^c \f_j)^n \rightarrow (\om+dd^c \f)^n$.
\end{pro}

\begin{proof}
We can assume without loss of generality that $\f \leq \f_j \leq 0$. Set
$$
\f_j^K:=\max (\f_j,-K) \text{ and } \f^K:=\max(\f,-K).
$$
It follows from the fundamental inequalities that the $(\chi,\om_j)$-energy
of the functions $\f_j^K$ is uniformly bounded independently of
both $j$ and $K$. 

Let $\nu$ be any cluster point of the bounded sequence of positive measures
$\nu_j:=(-\chi) \circ \f_j^K \, (\om_j+dd^c \f_j^K)^n$.
Observe that $(\om_j+dd^c \f_j^K)^n \rightarrow (\om+dd^c \f^K)^n$
as $j \rightarrow +\infty$, since $\f_j^K$ decreases towards $\f^K$ 
and these functions are uniformly bounded [BT 3]. It follows therefore
from the upper-semi-continuity of $\f^K$ that
$$
0 \leq (-\chi) \circ \f^K \, (\om+dd^c \f^K)^n \leq \nu.
$$
In particular $E_{\chi}(\f^K)$ is bounded from above
by $\nu(X)$, hence $\f \in\E$ by Corollary 2.4.

The convergence of $(\om_j+dd^c \f_j)^n$ towards $(\om+dd^c \f)^n$
follows again from the fact that
$(\om_j+dd^c \f_j^K)^n-(\om_j+dd^c \f_j)^n$ converges towards zero as 
$K \rightarrow +\infty$, uniformly with respect to $j$.
\end{proof}

We refer the reader to [EGZ] for an application to the construction
of K\"ahler-Einstein metrics on canonical/minimal models in the sense of Mori.

\vskip .2cm

Vincent Guedj \& Ahmed Zeriahi

Laboratoire Emile Picard

UMR 5580, Universit\'e Paul Sabatier

118 route de Narbonne

31062 TOULOUSE Cedex 04 (FRANCE)

guedj@picard.ups-tlse.fr

zeriahi@picard.ups-tlse.fr

\end{document}